 \documentstyle [12pt,twoside, epsf]{article}

\newcommand{\ZZ}{\mbox{$Z\!\!\! Z\!$}}	
\newcommand{\CC}{\mbox{$C\!\!\!\! I$}}	

\let\\=\cr

\def\Chi{\hbox{\raise0.5ex\hbox{$\chi$}}}

\newtheorem{th}{Theorem}
\newtheorem{lem}{Lemma}
\newtheorem{cor}{Corollary}
\newtheorem{prop}{Proposition}

\newtheorem{defn}{Definition}

\newtheorem{rem}{Remark}

\def\picill#1by#2(#3){\epsffile{#3}}
 
\begin{document} 
\pagestyle{myheadings}
\markboth{{\sc Kauffman \& Lambropoulou}}{{\sc Rational tangles and continued fractions}}

\title{On the classification of rational tangles}




\author{Louis H. Kauffman \\
Department of Mathematics, Statistics and
Computer Science\\
University of Illinois at Chicago\\
851 South Morgan St., Chicago IL 60607-7045, USA\\
kauffman@uic.edu\\
and\\
Sofia Lambropoulou \\
National Technical University of Athens\\
Department of Mathematics\\
Zografou campus, GR-157 80 Athens, Greece\\
sofia@math.ntua.gr
}

\date{}

\maketitle


\begin{abstract} 
\noindent In this paper we give two new combinatorial proofs of the classification of rational
tangles using the calculus of continued fractions. One proof uses the  classification of
alternating knots. The other proof uses  colorings of tangles. We also obtain an elementary
proof that alternating rational tangles have minimal number of crossings. Rational tangles form a
basis for the classification of knots and are of fundamental importance in the study of DNA
recombination. 
\end{abstract}


\section{Introduction}

A   {\it rational tangle}  is a proper embedding of two unoriented arcs  $\alpha_1, \alpha_2$  in a
$3$-ball $B^3,$ so that the four endpoints lie in the boundary of $B^3,$  and such that there
exists a homeomorphism of pairs: 

$$\overline{h}: (B^3, \alpha_1, \alpha_2) \longrightarrow (D^2\times I, \, \{x,y\}\times I) {\mbox
\ \ (a \ trivial \ tangle).}$$
  
\noindent This is equivalent to saying that rational tangles have specific representatives obtained
by applying a finite number of consecutive twists of neighbouring endpoints  starting from two 
unknotted and unlinked arcs (see Note 1 in Section 2). Such a pair of arcs comprise the $[0]$ or
$[\infty]$ tangles, depending on their position in the plane (Figures 1 and 2). We shall use this
characterizing  property of rational tangles as our definition (Definition 1 below).

\bigbreak

 We  are interested in tangles up to isotopy. Two rational tangles, $T, S$, in $B^3$ are {\it
isotopic}, denoted by
$T
\sim S$,  if there is an orientation-preserving
self-homeomorphism $h: (B^3, T) \longrightarrow (B^3, S)$ that is the identity map on the boundary.
Equivalently, $T, S$ are isotopic if and only if any two  diagrams of theirs (i.e. seeing the 
tangles as planar graphs) have identical configurations of their four endpoints on the boundary of
the projection disc, and they differ by a finite sequence of the well-known Reidemeister moves
\cite{Rd1}, which take place in the interior of the  disc.  Of course, each twisting operation
changes the isotopy class of  the tangle to which it is applied.  

\bigbreak
The rational tangles consist in a special class of
$2$-tangles, i.e. embeddings in a $3$-ball of two arcs and a finite number of circles.
The $2$-tangles  are  particularly interesting  because of the simple symmetry of their endpoints, 
which keeps the class closed under the tangle operations (see Figure 3  below). 
Moreover, the special symmetry of the endpoints of $2$-tangles allows for the following
{\it closing} operations, which yield two different knots or links: The {\it Numerator} of a 
$2$-tangle,
$T$, denoted by $N(T)$, which is obtained by joining with simple arcs the two upper endpoints and 
the two lower endpoints of $T,$ and the {\it Denominator} of a $2$-tangle, $T$, which is obtained
by joining with simple arcs each pair of the corresponding top and bottom endpoints of $T$,
and it shall be denoted by $D(T)$. Every knot or link can arise as the numerator closure of a
$2$-tangle.  The theory of general tangles has been introduced in 1967 by  John H. Conway \cite{C1}
in his work on enumerating and classifying knots. (In fact Conway had been thinking about tangles 
since he was a student in high school and he obtained his results as an undergraduate student in
college.)

\bigbreak
 The rational tangles give rise via numerator or denominator closure to a special class
of knots and links, the {\it rational knots} (also known as Viergeflechte, four-plats and
$2$-bridge knots).  These have one or two components, they are alternating and they are the
easiest knots and links to make (also for Nature, as DNA recombination suggests). The first twenty
five knots, except for
$8_5$, are rational. Furthermore all knots and links up to ten crossings are either rational or
are  obtained by inserting rational tangles into a few simple planar graphs, see \cite{C1}. The
2-fold branched covering spaces of $S^3$ along the rational knots  give rise to the lens spaces
$L(p,q)$  \cite{Sei}, \cite{ST}. Different rational tangles can
give the same rational knot when closed and this leads to  the subtle theory of the classification
of rational knots, see \cite{Sch}, \cite{Bu} and \cite{KL1}.  Finally,
rational knots and rational tangles figure prominently in the applications of knot theory to the
topology of DNA, see \cite{ES}, \cite{Su}.  Treatments  of various  aspects of rational tangles and
rational knots can be found in various places in the literature, see \cite{C1}, \cite{Sie}, \cite{BZ},
\cite{Kaw}, \cite{M}. See also 
\cite{BM} for a good discussion on classical relationships of rational tangles, covering spaces and
surgery. At the end of the paper we give a short history of rational knots and rational tangles.  

$$ \picill2.1inby2.35in(ratl1) $$

\begin{center}
{Figure 1 - A rational tangle in standard form } \end{center}
\vspace{3mm}

\noindent  A rational tangle  is associated in a canonical manner with a unique, reduced  rational
number or $\infty,$ called {\it the fraction} of the tangle. Rational tangles are classified by 
their fractions by means of the following theorem due to John H. Conway \cite{C1}:

\begin{th}[Conway, 1970]{ \ Two rational tangles are isotopic if and only if they have the 
same fraction. 
 } \end{th}

In  \cite{C1} Conway defined the fraction of a
rational tangle using its continued fraction form.  
He also defined a topological  invariant $F(R)$ for an arbitrary $2$-tangle $R$ using the
Alexander polynomial of the  knots  $N(R)$ and $D(R),$ namely as:  $F(R) = \frac{\Delta
(N(R))}{\Delta (D(R))}.$ He then observed that this evaluated at $-1$ coincides with the fraction
for rational tangles. The advantage of the second definition is that it is already a topological
invariant of the tangle. Proofs of Theorem~1 are given in
\cite{Mo}, \cite{BZ} p.196 and \cite{GK2}.  The first two proofs  used the second 
 definition of the fraction as an isotopy invariant of rational tangles. Then, for proving that
the fraction classifies the rational tangles, they invoked the classification of rational knots. 
 The proof by Goldman
and Kauffman \cite{GK2} is the first combinatorial proof of the classification of rational tangles. In
\cite{GK2}  the fraction of an unoriented $2$-tangle $S$ is defined via the bracket polynomial of the
unoriented knots $N(R)$ and $D(R),$  namely as: $F(R) = i\frac{<N(R)>(A)}{<D(R)>(A)},$  where
the indeterminate $A$ is specified to
$\sqrt{i}.$  There again the fraction is by definition an isotopy invariant of the tangles. 
The first definition of the fraction is more natural, in the sense that it is obtained directly from
the topological structure of the rational tangles. In order to prove Theorem 1 using this 
definition we need to  rely on a deep result in knot theory -- namely the
solution  of the {\it Tait Conjecture}
\cite{Ta}  concerning the classification of alternating knots that was given by Menasco and
Thistlethwaite \cite{MT} in 1993, and to adapt it to rational tangles. 
\bigbreak

 It is the main purpose of this paper to give this direct combinatorial proof of Theorem~1.  We
believe that our proof gives extra insight into the isotopies of rational tangles and the nature
of the theorem beyond the proof in \cite{GK2}.   The fraction is defined directly from the
algebraic combinatorial structure of the rational tangle by means of a continued fraction
expansion, and we have to show that it is an isotopy invariant. The topological invariance of the
fraction is proved via flyping.  
 We will show that the fraction is invariant under flyping
(Definition 2) and the transfer moves (see Figure 14), from which it follows that it is an isotopy
invariant of rational tangles. We will also show that two rational tangles with the same  fraction
are isotopic. These two facts imply Theorem~1.  

In the
course of this proof we will see and we will exploit the extraordinary interplay between   the
elementary number theory of continued fractions and the topological structure of rational tangles,
using their  characteristic properties:  The rational flypes (Definition 2) and equivalence of
flips (Definition 3).  The core of our proof is that rational tangles and continued fractions
have a similar canonical form, and the fact that rational tangles are alternating, for which we
believe we found the simplest possible proof. This implies the known result that the rational
knots are alternating.   We also give a second combinatorial proof of Theorem 1 
by defining in section 5, without using the Tait conjecture, the tangle fraction
via coloring. This paper serves as a basis for a sequel  paper
\cite{KL1}, where we give the first combinatorial proofs of Schubert's classification theorems for
unoriented and oriented rational knots  \cite{Sch}, using the results and the techniques developed
here. 

\bigbreak

 The paper is organized as follows. In Section 2 we introduce the operations on rational
tangles, we discuss the Tait conjecture for alternating knots and we prove a canonical form for 
rational tangles. In Section 3 we discuss some facts about continued fractions and we prove a key
result, a unique  canonical form.  In Section 4 we define the fraction of a rational tangle, we  
unravel  in full the analogy between continued fractions and rational tangles (analogy of
operations and calculus), and we give our proof of the classification of rational tangles. We also
prove the minimality of crossings for alternating rational tangles without necessarily resting on
the solution to the Tait conjecture. In Section 5 we give an alternate definition of the fraction
of a rational tangle via  integral coloring, as well as another combinatorial proof of Theorem 1,
without using the Tait conjecture. In Section 5 we use the structure of integral colorings of rational tangles
to prove for rational knots and links a special case of a conjecture of Kauffman and Harary \cite{HK} about colorings of
alternating links. Finally, in
Section 6 we reduce the number of operations that generate the rational tangles and we give a short history of rational knots and
rational tangles.  Throughout the paper by `tangle'  we will mean `tangle diagram' and by `knots' we will be
referring to both knots and links.

\section{The Canonical Form of Rational Tangles}

Clearly the simplest rational tangles are the
$[0]$, the
$[\infty]$, the $[+1]$ and the 
$[-1]$ tangles, whilst the next simplest ones are: 
 \begin{itemize}
\item[(i)] The  {\it integer tangles}, denoted by $[n],$ made of $n$ horizontal twists, $n \in
\ZZ,$ 
\vspace{-.05in}
\item[(ii)] The {\it vertical tangles},  denoted by $\frac{1}{[n]},$ made of $n$ vertical twists, 
$n
\in \ZZ.$
\end{itemize} 

$$ \picill5.3inby1.85in(ratl2) $$

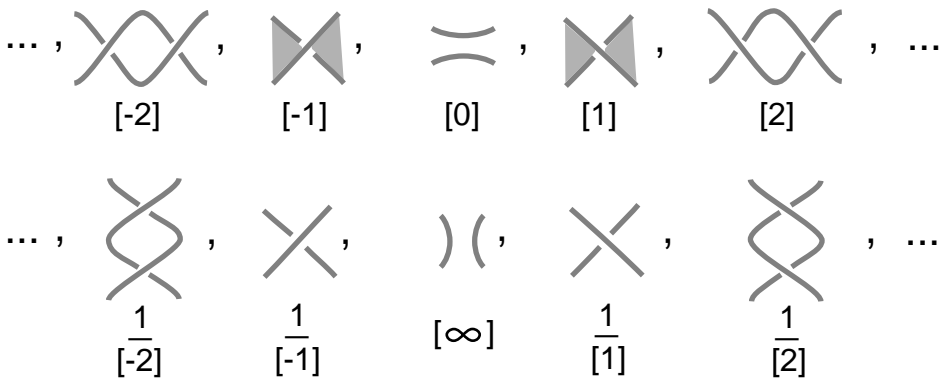
\begin{center} {Figure 2 - The elementary rational tangles and the shading rule} 
\end{center}
\vspace{3mm}

\noindent We note that the type of crossings of knots and tangles follow
the checkerboard  rule: Shade the regions of the tangle (knot) in two colors, starting from the
left (outside) to the right (inside) with grey, and so that adjacent regions have different
colors.  Crossings in the tangle are said to be of {\it positive type} if they are arranged with
respect to the shading  as exemplified in Figure 2 by the tangle $[+1],$ whilst crossings of the
reverse type are said to be of {\it negative type} and they are exemplified in Figure 2 by the
tangle $[-1].$ The reader should note that our crossing type  conventions are
the opposite of those of Conway in \cite{C1} and of those of Kawauchi in \cite{Kaw}. Our conventions
agree with those of  Ernst and Sumners in \cite{ES}, which also follow the standard conventions of 
biologists.

\smallbreak
 Rational tangles can be added, multiplied, rotated, mirror imaged and inverted. These are
well-defined (up to isotopy) operations in the
class of $2$-tangles, adequately described in Figure 3. In particular, {\it the
sum} of two $2$-tangles is denoted by `$+$' and {\it the product} by `$*$'. Notice that addition
and multiplication of tangles are not commutative. Also, they do not preserve the
class of rational tangles. For example, the tangle $\frac{1}{[3]} + \frac{1}{[2]}$ is not rational. 
 We point out that the numerator (denominator) closure of the sum (product) of two rational tangles is
still a rational knot, but the sum (product) of two rational tangles is a rational tangle if and only if
one of the two is an integer (a vertical) tangle. 

 {\it The mirror image} of a tangle $T,$  denoted $-T,$ is obtained from
$T$ by switching  all the crossings. E.g. $-[n] = [-n]$ and 
$-\frac{1}{[n]} = \frac{1}{[-n]}.$  Then we have $-(T+S)= (-T)+(-S)$ and $-(T*S)= (-T)*(-S).$ 
 Finally, {\it the rotation} of $T$, denoted 
$T^r$, is obtained by rotating $T$ counterclockwise by $90^0,$ whilst  {\it the inverse} of $T$,
denoted $T^i$,   is defined to be  ${-T}^r$. For example, ${[n]}^i = \frac{1}{[n]}$ and
${\frac{1}{[n]}}^i = {[n]}.$  Turning the tangle clockwise by $90^0$ is the cancelling operation of
our defined inversion, denoted $T^{-i}.$  In particular ${[0]}^r = {[0]}^i = [\infty]$ and ${[\infty]}^r
= {[\infty]}^i = [0].$   We have that $N(T) = D(T^r)$ and
$D(T) = N(T^r).$

 $$\vbox{\picill4.9inby1.9in(ratl3)  }$$

\begin{center}
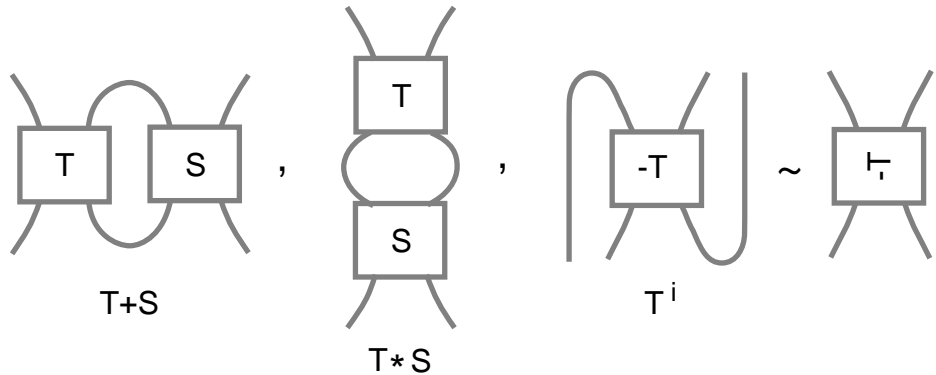
{Figure 3 - Addition, multiplication and inversion of $2$-tangles } 
\end{center}
\vspace{3mm} 

\noindent  Note  that $T^r$ and $T^i$ are in general not isotopic to
$T$. Also, it is in general not the case
that the inverse of the inverse of a $2$-tangle is isotopic to the original tangle, since
$(T^{i})^{i} = (T^{r})^{r}$ is the tangle obtained from $T$ by rotating it on its plane by $180^0$.
 For $2$-tangles the inversion is an order four operation. But, remarkably, for rational tangles the
inversion is an operation of order two, i.e. $T^{-i} \sim T^i$ and $T \sim (T^{i})^{i}$ (see Lemma 2).  
For this reason we shall denote the inverse of a rational tangle $S$ as
$\frac{1}{S}.$ This explains the notation for the vertical tangles. 
In particular we shall have $\frac{1}{[0]} = [\infty]$ and $\frac{1}{[\infty]} = [0].$

\begin{defn}{\rm \  A rational tangle is in {\it twist form} if it is created by consecutive
additions  and multiplications by the tangles $[\pm 1]$, starting  from
the tangle $[0]$ or the tangle $[\infty].$ (See Figure 4 for an example.) 
 } \end{defn}
Conversely, a rational tangle in twist form can be  brought to one of the  tangles $[0]$ or 
$[\infty]$  by a finite sequence of untwistings.  It follows that a rational tangle in twist form
can be obtained inductively from a  previously created rational tangle  by consecutive additions of
integer tangles and multiplications by vertical tangles, and it can be described by an algebraic
expression of the type:  
\[[s_k] +( \cdots + (\frac{1}{[r_3]} * ([s_1] + (\frac{1}{[r_1]} * [s_0] * \frac{1}{[r_2]}) + [s_2])
* \frac{1}{[r_4]}) + \cdots )+ [s_{k+1}],\] 
\noindent or of the type: 
\[\frac{1}{[r_k]} * ( \cdots * ([s_3] + (\frac{1}{[r_1]} * ([s_1] + \frac{1}{[r_0]} +
[s_2]) * \frac{1}{[r_2]}) + [s_4]) * \cdots ) * \frac{1}{[r_{k+1}]},\] 
according as we start building from the tangle
$[0]$ or  $[\infty],$ where all $s_i, r_i \in \ZZ.$ Note that some of the $s_i, r_i$ may be zero. 
By allowing $[s_k] + [s_{k+1}] = [0]$ and $[s_0] = [\infty]$ in the first expression, an algebraic
expression of the following type can subsume both cases. 

\[ T=[s_k] +( \cdots + (\frac{1}{[r_3]} * ([s_1] + (\frac{1}{[r_1]} * [s_0] * \frac{1}{[r_2]}) +
[s_2]) * \frac{1}{[r_4]}) + \cdots )+ [s_{k+1}],\] 

\noindent where $s_i, r_i \in \ZZ.$   For example,
the rational tangle of Figure 4 can be described as \  $( ([3] + ([1] * [3] * \frac{1}{[2]}) + [-4]) *
\frac{1}{[-4]}) + [2].$ With the above notation and for any $j \leq k$ we call {\it a truncation of
$T$} the result of
 untwisting $T$ for a while, i.e. a rational tangle of the type: 

\[ R = [s_j] +( \cdots + (\frac{1}{[r_3]} * ([s_1] + (\frac{1}{[r_1]} * [s_0] * \frac{1}{[r_2]}) +
[s_2]) * \frac{1}{[r_4]}) + \cdots )+ [s_{j+1}].\]

$$ \picill3.55inby3.6in(ratl4) $$

\begin{center}
{Figure 4 - A rational tangle in twist form}
\end{center}
\vspace{3mm}

\noindent {\bf Note 1} \ To see the equivalence of Definition 1 with the definition of a rational
tangle given in the introduction let
$S^2$ denote the two-dimensional sphere, which is the boundary of the $3$-ball, $B^3,$  and let $p$
denote four specified points  in
$S^2.$ Let further \ $h: (S^2,p) \longrightarrow (S^2,p)$ \ be a self-homeomorphism of
$S^2$  with the four points. This extends to a self-homeomorphism
$\overline{h}$ of the $3$-ball $B^3$ (see \cite{R}, page 10).  Further, let $a$ denote the
two straight arcs $\{x,y\}\times I$  joining pairs of the fours point of the boundary of
$B^3.$ Consider now $\overline{h}(a).$ We call this the tangle induced by $h.$ We note that,
up  to isotopy,  $h$ is a composition of braidings of pairs of
points in $S^2$ (see \cite{PS}, pages 61 to 65).  Each such braiding induces a twist in
the corresponding tangle. So, if $h$ is a composition of braidings of pairs of points,
then the extension $\overline{h}$ is a composition of twists of neighbouring end arcs. 
Thus $\overline{h}(a)$ is a rational tangle and every rational tangle can be obtained this
way.  
\bigbreak

 We define now an isotopy move for
rational tangles that plays a crucial role in the whole
 theory that follows. 

\begin{defn}{\rm \  A {\it flype} is an isotopy of a $2$-tangle/a knot applied on a $2$-subtangle
of the form $[\pm 1]+t$ or $[\pm 1]*t$ as illustrated in Figure 5. A flype fixes the endpoints of the
subtangle on which it is applied.  A flype shall be called {\it rational} if the
$2$-subtangle on which it acts is rational. }
\end{defn}

$$ \picill3.7inby2.2in(ratl5) $$

\begin{center}
{Figure 5 - The flype moves} 
\end{center}
\vspace{3mm} 

\noindent  A tangle is said to be {\it alternating} if the crossings 
alternate from under to over as we go along any component or arc of the weave.  Similarly, a knot is
{\it alternating} if it possesses an alternating diagram.
 Notice that, according to the checkerboard shading, the only way the weave alternates is
if any two adjacent crossings are of the same type, and this propagates to the whole diagram. Thus,
{\it a tangle or a knot diagram with all crossings of the same type is alternating}, and
this characterizes alternating tangles and knot diagrams. It is important to note that flypes preserve
the alternating structure. Moreover, flypes are the only isotopy moves needed in the statement of the
celebrated Tait Conjecture for alternating knots. This was P.G. Tait's  working assumption in 1877
(see
\cite{Ta}) and was proved   by W.~Menasco and M.~Thistlethwaite \cite{MT} in 1993. 
\bigbreak
\noindent {\bf The Tait Conjecture for Knots}. Two alternating knots
 are isotopic if and only if any two corresponding diagrams on $S^2$ are related by a finite
sequence of flypes.  
\bigbreak

\noindent For rational tangles flypes are of very specific types, as the lemma below shows.
\begin{lem}{ \ Let  $T$ be a  rational tangle in twist form. Then
 \begin{itemize}
\item[(i)] \ $T$ does not contain any non-rational $2$-subtangles.  
\vspace{-.1in}
\item[(ii)] \ Every $2$-subtangle of $T$ is a truncation of $T$.
\end{itemize} 
 } \end{lem}

\noindent {\em Proof.} \  By induction. Notice that both statements are true for the tangles $[0],
[\infty]$ and $[\pm1].$ Assume they are true for all rational tangles with less than $n$ crossings,
and let $T$ be a rational tangle in twist form with $n$ crossings. By Definition 1 the tangle $T$
will contain an outmost crossing, i.e. $T = T'+ [\pm1]$ or $T = [\pm1] +T'$ or $T = T'* [\pm1]$
or $T = [\pm1] * T'.$ 

For proving $(i)$ we proceed as follows. Let  $U$ be a $2$-subtangle of
$T.$ Then $U$ either contains the outmost crossing of $T$ or not. If $U$ does not contain the
crossing, then by removing it we have $U$ as a $2$-subtangle of the tangle $T'.$ But $T'$ has $n-1$
crossings, and by induction hypothesis $U$ is rational. If $U$ does  contain the outmost crossing,
then by removing it we also remove it from $U,$ and so we obtain  a $2$-subtangle $U'$ of the new
tangle $T'.$ But $U$ is rational if and only if $U'$ is rational, and $U'$ has to be rational by
induction hypothesis. 

 For proving $(ii)$ let $U$ be a $2$-subtangle of $T.$  By $(i)$ $U$ has to be
rational and, arguing as in $(i),$ $U$ either contains the outmost crossing of
$T$ or not. If not, then by removing the crossing we have $U$ as a $2$-subtangle
of the tangle $T',$ and by induction hypothesis $U$ is a truncation of $T',$ and thus also of $T.$ 
If $U$ does  contain the outmost crossing, then by removing it we 
obtain  a $2$-subtangle $U'$ of the new tangle $T',$ and by induction hypothesis $U'$ is a
truncation of  $T'.$ Then $U'$  is also a truncation of $T,$ and thus so is $U.$ 
 $\hfill \Box $

\begin{cor}{ \ All flypes of a rational tangle $T$ are rational.  }
\end{cor}   

\begin{defn}{\rm \  A {\it flip} is a rotation in space of a $2$-tangle by $180^0$. We say that
$T^{hflip}$ is the {\it horizontal flip} of the $2$-tangle $T$ if $T^{hflip}$ is obtained from $T$
by a $180^0 - $ rotation around a horizontal axis on the plane of $T$, and $T^{vflip}$ is the {\it
vertical flip} of the tangle $T$ if $T^{vflip}$ is obtained from $T$ by a $180^0 - $ rotation around a
vertical axis on the plane of $T$, see Figure 6 for illustrations.  
} \end{defn}  

$$ \picill3.63inby2.35in(ratl6) $$

\begin{center}
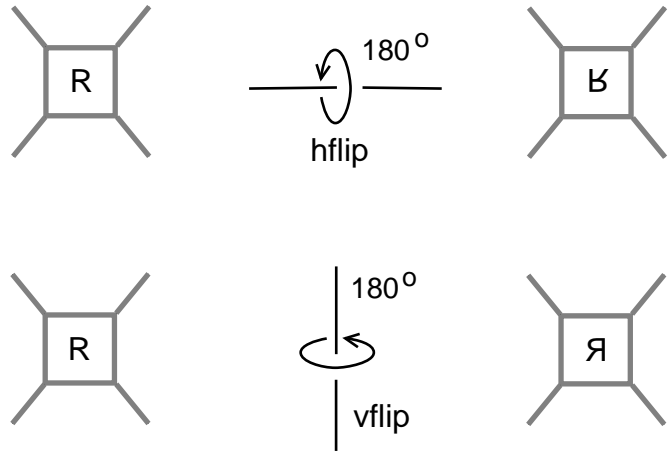
{Figure 6 - The horizontal and the vertical flip} 
\end{center}
\vspace{3mm}

\noindent In view of the above definitions, a flype on a $2$-subtangle $t$ can be
described by  one of the isotopy identities: 
\[ [\pm 1] + t \sim t^{hflip} +[\pm 1] \mbox{ \ \ \ \ or \ \ \ \ } [\pm 1] * t \sim t^{vflip} * 
[\pm 1].\]

\noindent Now we come to a remarkable property of rational tangles. Note that a flip switches the
endpoints of the tangle and, in general, a flipped tangle is not isotopic to the original one.
But this is the case for rational tangles, as the lemma below shows.

\begin{lem} {\bf (Flipping Lemma)}{ \ If \, $T$ is rational, then: 
\[ (i) \ T \sim T^{hflip}, \ \ \ \ (ii) \ T \sim T^{vflip} \ \ \ \ \mbox{and} \ \ \ \ (iii) \ T
\sim (T^{i})^{i} = (T^{r})^{r}.\]
 } \end{lem}

\noindent {\em Proof.} We prove (i) and (ii) by induction. Note that both statements are true for
the tangles $[0], \ [\infty], \ [\pm 1],$ and assume they are true for any rational tangle,
$R$ say, with $n$ crossings, i.e. $R \sim R^{hflip}$ and $R \sim R^{vflip}.$  We will show that
then the statements hold also for the tangles $F~=~R~+~[\pm 1],$ \ $F'~=~[\pm 1]~+~R,$ \ 
$L~=~R~*~[\pm 1],$ \ and $L'~=~[\pm 1]~*~R.$ 
Then, by Definition~1 and by Note 1, the statements shall be true for any rational tangle. Indeed, 
for 
$F^{hflip}$ and $L^{hflip}$ we have:

$$ \picill15cmby2.5in(ratl7) $$

\begin{center}
{Figure 7 - The proof of Lemma 2} 
\end{center}
\vspace{3mm} 

\noindent With the same arguments we show that $F^{vflip} \sim F$ and $L^{vflip} \sim L.$ For the 
tangles $F'$ and $L'$ the proofs are completely analogous. Finally,
statement (iii) follows from (i) and (ii), since  $(T^{i})^{i}=(T^{r})^{r}=(T^{hflip})^{vflip}.$
 $\hfill \Box $

\begin{rem}{\rm \   As a consequence of Lemma 2, addition of $[\pm 1]$ and multiplication by  $[\pm 1]$
are commutative, so a rational flype is described by
\[ [\pm 1] + t \sim t +[\pm 1] \mbox{ \ \ \ \ or \ \ \ \ } [\pm 1] * t \sim t * [\pm 1]. \]
\noindent In general for any $m, n \in \ZZ$ we have the following isotopy
identities: 
\[ [m]+T+[n] \, \sim \, T+[m+n], \ \ \ \ \frac{1}{[m]}*T*\frac{1}{[n]} \, \sim \, T*\frac{1}{[m+n]}.
\] 

\noindent In view of Lemma 2, another way to define a rational flype is by  one of the following isotopy
identities: 
\[ [\pm 1] + t \sim ([\pm 1] + t)^{vflip} \mbox{ \ \ \ \ or \ \ \ \ } [\pm 1] * t \sim ([\pm 1] *
t)^{hflip}.\]
} \end{rem}

  Lemma 2(iii) says that inversion is an operation of order 2 for rational
tangles. Thus, if $T$  rational then $T^{i} \sim T^{-i},$ so we can rotate the mirror image of $T$ 
 by $90^0$ either counterclockwise 
 or clockwise to obtain $T^i$.  Thus, for a rational tangle $T$ its inverse shall be denoted
by $\frac{1}{T}$ or $T^{-1}$.  With this notation we have  $\frac{1}{\frac{1}{T}} =
T$ and $T^{r} = \frac{1}{-T} = -\frac{1}{T}.$

\begin{defn}{\rm \ A rational tangle is said to be in {\it standard form} if it is 
created by consecutive additions of the tangles $[\pm 1]$  {\it only on the right} (or only on the left)
and multiplications  by the tangles $[\pm 1]$ {\it only at the bottom} (or only at the top), starting 
from the tangle
$[0]$ or $[\infty].$    }
\end{defn}

\noindent  Thus, a rational tangle in standard form can be obtained inductively from a 
previously created rational tangle, $T$ say, either by adding an integer tangle  on the right:
$T \rightarrow T+[\pm k],$ or by  multiplying  by a vertical tangle  at the bottom: 
$T~\rightarrow~T~*~\frac{1}{[\pm k]}$, starting from $[0]$ or $[\infty],$ see Figure 8. 

$$ \picill3.8inby2.65in(ratl8) 
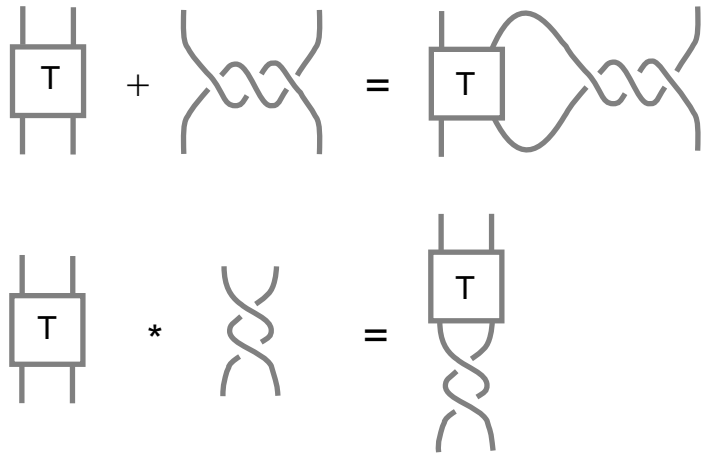
$$

\begin{center}
{Figure 8 - Creating rational tangles in standard form} 
\end{center}
\vspace{3mm} 

\noindent  Figure 1 
illustrates the tangle  $(([3]*\frac{1}{[-2]})+[2])$ in standard form. Hence, a rational tangle in
standard form has an algebraic expression of the type:  
\[ ( ( ([a_n] * \frac{1}{[a_{n-1}]}) + [a_{n-2}]) * \cdots * \frac{1}{[a_2]} ) + [a_1], \ \ 
 \mbox{for} \ a_2, \ldots, a_{n-1} \in \ZZ - \{0\}, \]

\noindent where $[a_1]$ could be $[0]$ and $[a_n]$ could be $[\infty]$ (see also Remark 2 below). 
The $a_i$'s are integers denoting numbers of twists with their types. Note that the tangle begins to
twist from the tangle  $[a_n]$ and it untwists from the tangle $[a_1].$

\noindent  Figure 9 illustrates two equivalent (by the Flipping Lemma) ways of  representing an
abstract rational tangle in standard form: The  {\it standard representation} of a
rational tangle. In either illustration the rational tangle begins  to twist from the
tangle $[a_n]$ ($[a_5]$ in Figure 9), and it untwists from the tangle $[a_1].$ 
Note that the tangle in Figure 9 has an odd number of sets of twists ($n=5$) and this causes 
$[a_1]$ to be horizontal.  If $n$ is even and $[a_n]$ is horizontal then $[a_1]$ has to be
vertical.

$$ \picill4.7inby2.2in(ratl9) $$

\begin{center}
{Figure 9 - The standard representations } 
\end{center}
\vspace{3mm} 

Another way of representing an abstract rational tangle in standard form is the {\it
$3$-strand-braid representation}, illustrated in Figure 10, which is more useful for studying
rational knots.  For an example see Figure 11. As Figure 10 shows, the $3$-strand-braid
representation is actually a compressed version of the standard representation, so the two
representations are equivalent. The upper row of crossings of the $3$-strand-braid  representation
corresponds to the horizontal crossings of the standard
representation and the  lower row to the
vertical ones, as it is easy to see by a planar rotation. Note that, even though
the type of crossings does not change by this planar rotation, we need to draw the mirror
images of the even terms, since when we rotate them to the vertical position we obtain
crossings of the opposite type in the local tangles. In order to bear in mind this 
change of the local signs we put on the geometric picture the minuses on the even terms.

$$ \picill5.5inby2.1in(ratl10) $$

\begin{center}
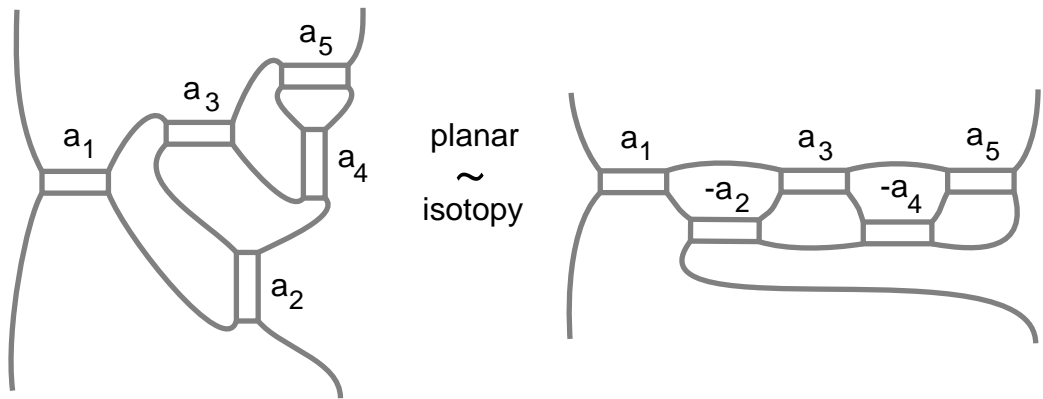
{Figure 10 - The standard and the 3-strand-braid representation } 
\end{center}
\vspace{3mm}

\begin{rem}{\rm \  When
we start creating a rational tangle, the very first crossing can be equally seen as
a horizontal or as a vertical one. Thus, we may always assume that we start twisting from the
$[0]$-tangle. Moreover, because of the same ambiguity, we may always assume that {\it the index
$n$} in the above notation {\it is always odd}. This is illustrated in Figure 11.
 }\end{rem}

$$ \picill5.3inby1.1in(ratl11) $$

\begin{center}
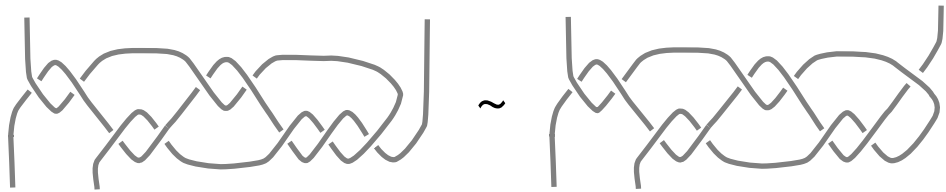
{Figure 11 - The ambiguity of the first crossing} 
\end{center}
\vspace{3mm} 

\noindent  From the above one may associate
to a rational tangle $T$ a vector of integers $(a_1, a_2, \ldots, a_n)$. The first entry denotes
the place where $T$ starts unravelling and the last entry is where it  begins to twist.
 For example the rational tangle of Figure 1 is associated to the induced vector $(2,
-2, 3),$  while the  tangle of Figure 4 corresponds after a sequence of flypes to
the vector $(2, -4, -1, 3, 3).$   For the rational tangle $T$  this
vector is {\it unique}, up to breaking the entry
$a_n$ by a unit, according to Remark 2. I.e. 
$(a_1, a_2, \ldots, a_n) = (a_1, a_2, \ldots, a_n -1, 1),$ if $a_n >0,$ and 
$(a_1, a_2, \ldots, a_n) = (a_1, a_2, \ldots, a_n +1, -1),$ if $a_n <0.$ (From the above  $n$ may 
be
 assumed to be odd.) As we shall soon see, if  $T$  changes by an isotopy the induced
associated vector is not the same.  

\bigbreak

\noindent  The following lemma shows that the standard form is generic for rational tangles.

\begin{lem}{ \ Every rational tangle can be brought via isotopy to standard form.
 } \end{lem}

\noindent {\em Proof.} Let $T$ be a rational tangle in twist form.   Starting from the outmost
crossings of $T$ and using horizontal and vertical rational flypes we bring, by induction, all
horizontal and all vertical twists to the right and to the bottom  applying the isotopy identities for
rational flypes given in Remark 1. This process yields that the tangle 
 
\[T = [s_k] +( \cdots + (\frac{1}{[r_3]} * ([s_1] + (\frac{1}{[r_1]} * [s_0] * \frac{1}{[r_2]}) +
[s_2]) * \frac{1}{[r_4]}) + \cdots )+ [s_{k+1}],\] 

\noindent  gets transformed isotopically to the tangle in standard form:

\[ (( ( ([s_0] * \frac{1}{[r_1 + r_2]}) + [s_1 + s_2]) *  \frac{1}{[r_3 + r_4]}) + \cdots ) + 
 [s_k + s_{k+1}]. \]
$\hfill\Box$
\bigbreak

\noindent For example, the tangle in Figure 4 is isotopic to the tangle $ ( ( ([3] * \frac{1}{[3]}) +
[-1]) *
\frac{1}{[-4]}) + [2]$ in standard form.

\begin{rem}{\rm \  It follows from Definition 4 and Lemma 3  that  the whole class of rational
tangles can be  generated inductively by the two simple algebraic operations below starting from
the tangles $[0]$ or $[\infty]$, where $T$ is any previously created rational tangle.   
  
\begin{enumerate}
\item {\it Right addition of $[+1]$ or $[-1]$:} \ \ \ \ \ \ \ \ \ \ \ $T \longrightarrow  T+[\pm
1].$ 
\item {\it Bottom multiplication by $[+1]$ or $[-1]$:} \ $T \longrightarrow T*[\pm 1].$ 
\end{enumerate}
} \end{rem}

\begin{defn}{\rm \ A {\it continued fraction in integer tangles} is an algebraic description of a
rational tangle via a continued fraction  built from the  tangles $[a_1],$  $[a_2],$ 
$\ldots, [a_n]$  with all numerators equal to~$1$, namely an expression of
the type:

\[ T =[[a_1],[a_2],\ldots ,[a_n]] \, := \, [a_1]+ \frac{1}{[a_2]+\cdots + \frac{1}{[a_{n-1}]
+\frac{1}{[a_n]}}}
\]

\noindent for $a_2, \ldots, a_n \in \ZZ - \{0\}$ and $n$ even or odd. We allow that the term
$a_1$  may be zero, and in this case the tangle $[0]$ may be omitted.  A rational tangle described
via a continued fraction in integer tangles is said to be in {\it continued fraction form}. The
{\it length} of the continued fraction is arbitrary -- here illustrated at length
$n$ whether the first summand is the tangle $[0]$ or not.  }
\end{defn}

\begin{lem}{ \ Every rational tangle $T$ satisfies the following isotopic equations:
\[ T* \frac{1}{[n]} = \frac{1}{[n] + \frac{1}{T}} \mbox{ \ \ \ \ and \ \ \ \ } 
\frac{1}{[n]} * T = \frac{1}{\frac{1}{T} + [n]}. \]
 } \end{lem}

\noindent {\em Proof.}  Figure 12 illustrates the proof of the first equation. Here `L.2' stands
for `Lemma 2'.  The second one is similar. That the two equations are indeed isotopic follows
from the proof of Lemma 3. $\hfill \Box $ 

$$ \picill15cmby1.75in(ratl12) $$

\begin{center}
{Figure 12 - The proof of Lemma 4} 
\end{center}
\vspace{3mm}

\begin{rem}{\rm \  It follows now from Remark 3 and Lemma 4 that the two simple
algebraic operations below generate inductively the whole class of rational tangles starting from
the tangle $[0]$, where $T$ is any previously created rational tangle.  
  
\begin{enumerate}

\item {\it Right addition of $[+1]$ or $[-1]$:} \ $T \longrightarrow T+[\pm 1].$ 
\item {\it Inversion of rational tangles:}  \ \ \ $T \longrightarrow \frac{1}{T} = T^{-1}.$ 
\end{enumerate}

\noindent  It is easy to see
that the second operation can be replaced by the operation: 

\begin{itemize} 
\item[2$'$.] {\it Rotation of rational tangles:} \ \  \ $T \longrightarrow T^{r} = -\frac{1}{T}.$
\end{itemize}
 } \end{rem}
 
\noindent In Section 6 we sharpen this even more by showing that the class of rational tangles is
generated inductively from the tangle $[0]$ by addition of $[+1]$  and rotation. We are now in a
position to prove the following:

\begin{prop}{ \  Every rational tangle can be written in continued fraction form.
} \end{prop}  

\noindent {\em Proof.}  By Lemma 3, a rational tangle may be assumed to be in standard form and so
by repeated applications of Lemma 4  we obtain the corresponding continued
fraction form:

\[ ( ( ([a_n] * \frac{1}{[a_{n-1}]}) + [a_{n-2}]) * \cdots * \frac{1}{[a_2]} ) + [a_1] \ 
\longrightarrow \ [a_1]+ \frac{1}{[a_2]+\cdots + \frac{1}{[a_{n-1}] +\frac{1}{[a_n]}}}. \] 
$\hfill \Box $

\noindent Thus the continued fraction form and the standard form of a rational tangle are
equivalent and the above correspondence shows that it is straightforward to write out the one from
the other.  For example, the tangle of Figure 1 can be written as 
$[2]+\frac{1}{[-2]+\frac{1}{[3]}},$  the one of Figure 4 as  
$[[2], [-4], [-1], [3], [3]],$  whilst the illustrations of  Figures 9 and 10 depict an abstract 
rational tangle
$[[a_1], [a_2], [a_3], [a_4], [a_5]].$  The following statements, now, about the continued fraction
form of rational tangles are straightforward.

\begin{lem}{ \ Let $T=[[a_1], [a_2],\ldots ,[a_n]]$ be a rational tangle in continued fraction form.
Then
 
\vspace{.15in}

\noindent $\begin{array}{lrcl} 
 
 1. & T + [\pm 1] &  = & [[a_1 \pm 1], [a_2],\ldots ,[a_n]],   \\    [1.8mm] 

 2. & \frac{1}{T} &  = & [[0], [a_1], [a_2],\ldots ,[a_n]],   \\  [1.8mm] 

3. & -T &  =  & [[-a_1], [-a_2],\ldots ,[-a_n]], \\  [1.8mm] 

4. & \mbox{If} \ R &  =  & [[a_{i+1}],\ldots ,[a_n]], \ 
\mbox{then we write} \ T=[[a_1], \ldots ,[a_i], R],  \\  [1.8mm] 

\end{array}$

\noindent \, $5.$ \ \ \ \ \ \ If \ $a_i \ = \ b_i + c_i$ \ and \ $S = [[c_i], [a_{i+1}],\ldots,
[a_n]],$ \ then \vspace{.07in}

$ \ T     =  [[a_1], \ldots, [a_{i-1}], [b_i] + S] =   
[[a_1],\ldots, [a_{i-1}], [b_i], [0], [c_i], [a_{i+1}], \ldots, [a_n]].$
 } \end{lem}

Recall that  a rational tangle $[[a_1], [a_2],\ldots ,[a_n]]$ is 
alternating if {\it the $a_i$'s are all positive or all  negative}.

\begin{defn}{\rm \ A rational tangle $T=[[\beta_1], [\beta_2], \ldots, [\beta_m]]$ is in {\it
canonical form} if $T$ is alternating and $m$ is odd. Moreover, $T$ shall be called {\it
positive} or {\it negative} according to the sign of its terms. }
\end{defn}

\noindent  We note that if $T$ is alternating and $m$ even, then we can bring $T$ to canonical
form by breaking $[\beta_m]$ to $[sign(\beta_m)\cdot (|\beta_m| - 1)] + [sign(\beta_m)\cdot 1],$ 
by Remark 2, and thus,
$[[\beta_1], [\beta_2], \ldots, [\beta_m]]$ to $[[\beta_1], [\beta_2],\ldots ,[sign(\beta_m)\cdot
(|\beta_m| - 1)], [sign(\beta_m)\cdot 1]].$ Proposition 2 below is a key property of rational 
tangles.

\begin{prop}{ \  Every  rational tangle can be isotoped to canonical form.
} \end{prop}  

\noindent {\em Proof.}  Let $T$ be a rational tangle. By Proposition 1, $T$ may be assumed to be 
in continued fraction form, say $T= [[a_1], [a_2], \ldots, [a_n]].$  We will show that
$T\sim[[\beta_1], [\beta_2], \ldots, [\beta_m]],$ where all $\beta_i$'s are positive or all
negative. If
$T$ is non-alternating then the
$a_j$'s are not all of the same sign. Let  $a_{i-1}, a_i$  be the first pair of adjacent $a_j$'s 
of opposite sign, and let 
$a_{i-1} > 0$. Then a configuration of the following type, as illustrated in Figure 13
below, must occur for $i$ odd  or a similar one for $i$ even.  
 
$$ \picill11.5cmby3.35in(ratl13a) $$

\begin{center}
{Figure 13 - A non-alternating configuration} 
\end{center}
\vspace{3mm} 

\noindent  If $a_{i-1} < 0$ then similar 
configurations will occur, but with the signs of $a_1, \ldots, a_i$ switched. We remind  that the 
signs of $a_{i+1}, \ldots,  a_n$ are irrelevant, and we note that the subtangles $t$ and
$s$ are rational and in continued fraction form. Now, inside $s$ the arc connecting the two 
crossings of opposite signs can be isotoped in both types of configurations to yield a simpler
rational tangle $s'$ isotopic to $s.$ See Figure 14 for $i$ odd and for $i$ even respectively. 
Such an isotopy move shall be called  {\it a transfer move}. Since
$s$ is a rational tangle in continued fraction form, the upper left arc of $s$ joins directly to
the subtangle
$a_n,$ and thus it meets no other arcs of the diagram.  Hence, after the transfer move the
subtangle $s'$ has  one fewer crossing than $s$ so we can apply induction.  

$$ \picill15cmby1.35in(ratl14a) $$

\noindent or

$$ \picill15cmby1.15in(ratl14b) $$

\begin{center}
{Figure 14 - The transfer moves} 
\end{center}
\vspace{3mm} 

\noindent The above isotopies are reflected in the following tangle identities for the cases $i$ 
odd and
$i$ even respectively. There are similar identities for switched crossings. 

\[ s = (t+[-1]) * [+1] \stackrel{L.4}{=} \frac{1}{
 [+1] + \frac{1}{[-1]+t}}  \sim - \frac{1}{t} + [+1] = s', \mbox {\em \ if $i$ odd, and} \]
\[ s = (t* [-1]) +[+1] \stackrel{L.4}{=} [+1] + 
 \frac{1}{[-1] + \frac{1}{t}} \sim - \frac{1}{t} * [+1] \stackrel{L.4}{=} \frac{1}{[+1] - t}
= s', \mbox {\em \ if $i$ even.} \]

\noindent In terms of tangle continued fractions the above can be
expressed as follows: 
\begin{itemize}
\item[{\em If $i$ odd:}] We have from Figure 13 that $t = [[a_i+1], [a_{i+1}], \ldots, [a_n]],$ 

$s = [[0], [+1], [-1] + t] = [[0], [+1], [a_i], \ldots, [a_n]]$ and, from Figure 14, that 

$s' = [[+1], -t] = [[+1], -[a_i+1], -[a_{i+1}], \ldots, -[a_n]].$ 

And so,

$T = [[a_1], \ldots, [(a_{i-1} - 1) +1], \ldots,[a_n]] = [[a_1], \ldots, [a_{i-2}], [a_{i-1} - 1]
+ \frac{1}{s}]$

$  \ \ \ = [[a_1], \ldots, [a_{i-2}], [(a_{i-1}-1)+(+1)], [-1]+t]$

$ \ \stackrel{L.5(5)}{=} [[a_1], \ldots, [a_{i-2}], [a_{i-1}-1], [0], [+1], [-1]+t]
\Longleftrightarrow$

$ T = [[a_1], \ldots, [a_{i-2}], [a_{i-1}-1], [0], [+1], [a_i], \ldots, [a_n]],$

which gets isotopically transformed to

$T' = [[a_1], \ldots, [a_{i-2}], [a_{i-1} - 1] + \frac{1}{s'}] = [[a_1], \ldots, [a_{i-2}],
[a_{i-1} - 1], s']$

$ \ \ \ \ = [[a_1], \ldots, [a_{i-2}], [a_{i-1}-1], [+1], -t] \Longleftrightarrow$

$ T' = [[a_1], \ldots, [a_{i-2}], [a_{i-1}-1], [+1], -[a_i+1], -[a_{i+1}], \ldots,
-[a_n]].$

\item[{\em If $i$ even:}] Here we have $t = [[0], [a_i+1], [a_{i+1}], \ldots, [a_n]],$

$s = [[+1], [-1]*t] \stackrel{L.4}{=} [[+1], [-1], t] \stackrel{L.5(5)}{=} [[+1], [a_i],
 \ldots, [a_n]]$ and

$s' \stackrel{L.4}{=} [[0], [+1]+(-t)] = [[0], [+1], -[a_i+1], -[a_{i+1}], \ldots, -[a_n]].$ 

And so,

$T = [[a_1], [a_2], \ldots, [a_n]] = [[a_1], \ldots, [a_{i-2}], [a_{i-1} - 1] + s]$

$ \ \ \ = [[a_1], \ldots, [a_{i-2}], [(a_{i-1}-1)+(+1)], [-1]+\frac{1}{t}]$  

$ \ \stackrel{L.5(5)}{=} [[a_1], \ldots, [a_{i-2}], [a_{i-1}-1], [0], [+1], [-1]+\frac{1}{t}]
\Longleftrightarrow$ 

$ T = [[a_1], \ldots, [a_{i-2}], [a_{i-1}-1], [0], [+1], [a_i], \ldots, [a_n]],$

which gets isotopically transformed to

$T' = [[a_1], \ldots, [a_{i-2}], [a_{i-1} - 1] + s']  = [[a_1], \ldots, [a_{i-2}], [a_{i-1} - 1], 
\frac{1}{s'}] $

$ \ \ \ \ = [[a_1], \ldots, [a_{i-2}], [a_{i-1}-1], [+1] + (-t)] \Longleftrightarrow $

$ T' = [[a_1], \ldots, [a_{i-2}], [a_{i-1}-1], [+1], -[a_i+1], -[a_{i+1}], \ldots, -[a_n]].$
\end{itemize} 

\noindent Notice that the breaking of $T$ as well as the final tangle
$T'$ are the same in either case. Note also, that the total number of crossings in $T'$ is indeed
reduced by one. For the cases of the same configurations, but with the signs of $a_1, \ldots, a_i$
switched we have completely analogous formulae. Thus, by induction $T$
is isotopic to an alternating rational tangle $[[\beta_1], [\beta_2], \ldots, [\beta_m]],$  where
$m$ is odd by the discussion before the proposition.  

Finally observe that, if the above isotopy
involves the integer tangle $[a_1],$ the transfer move will not be needed again in the same
region. Thus, in principle, the sign of $a_1$ or of $a_2,$ if $a_1=0,$ dominates the type of
crossings in the alternating weave.  There is one exception to this rule, namely when the tangle
begins with an alteration of $[+1]$ and $[-1]$ tangles. More precisely, if $T = [[+1], [-1], t],$
then the sign of $T$ is opposite to the sign of $t.$ If $T = [[+1], [-1], [+1], [-1], t],$ then
the sign of
$T$ is same as the sign of $t$, and if $T = [[+1], [-1], [+1], [-1], [+1], [-1], t],$ then 
$T = t.$  There are analogous considerations for alterations of $[-1]$ and  $[+1]$. The proof is 
now completed.  $\hfill \Box$ 

\bigbreak \noindent  The alternating nature of the  rational  tangles will be  very useful to us 
in classifying rational knots in \cite{KL1}.  It is easy to see that the closure of an alternating
rational tangle is an alternating knot. Thus we have 
\begin{cor}{ \  Rational knots are alternating, since they possess a diagram that is the closure 
of an alternating rational tangle. }
\end{cor}

\section{Some facts about Continued Fractions}

It is clear that every rational number can be written as a continued fractions with all numerators
equal to $1,$  namely as an arithmetic expression of the type: 

\[ [a_1, a_2, \ldots, a_n] := a_1+ \frac{1}{a_2+\cdots + \frac{1}{a_{n-1} +\frac{1}{a_n}}}  \]

\noindent for $ a_1 \in \ZZ, \ a_2, \ldots, a_n \in \ZZ - \{0\}$ and $n$ even or odd. As in the
case of rational tangles we allow that the term $a_1$  may be zero.  In the case of the  subject 
at hand we shall only consider this kind of continued fractions.  The subject of continued 
fractions is of perennial interest to mathematicians, see for example \cite{Ko}, \cite{O},
\cite{W}, \cite{Kh}.  The {\it length} of the continued fraction is the number $n$ whether $a_1$
is zero or not.  Note that if  $a_1 \neq 0$ ($a_1 = 0$), then the absolute value of the continued
fraction is greater (smaller) than one. Clearly, the two simple algebraic operations {\it
addition of $+1$ or $-1$} and {\it inversion}  generate inductively the whole class of continued
fractions starting from zero.
\smallbreak
\noindent In this section we prove a well-known canonical form for continued fractions. The 
algorithm we develop works in parallel with the algorithm for the canonical form of rational
tangles in the previous section.   The following statements about continued
fractions are really straightforward (compare with Lemma 5).

\begin{lem}{ \ Let $\frac{p}{q}$ be any rational number. Then
\vspace{.15in}

\noindent \ $1.  \ \ \mbox{there are} \ a_1 \in \ZZ, \ a_2, \ldots, a_n \in \ZZ - \{0\} \
\mbox{such that} \ \frac{p}{q} = [a_1, a_2, \ldots, a_n],  $

\noindent $\begin{array}{lrcl} 
 
 2.& \frac{p}{q} \pm 1 &  = & [a_1 \pm 1, a_2,\ldots ,a_n],   \\    [1.8mm] 

 3. & \frac{q}{p} &  = & [0, a_1, a_2, \ldots, a_n],   \\  [1.8mm] 

4. & -\frac{p}{q} &  =  & [-a_1, -a_2, \ldots, -a_n], \\  [1.8mm] 

5. & \mbox{If} \ \frac{r}{d} &  =  & [a_{i+1},\ldots ,a_n], \ 
\mbox{then we write} \ \frac{p}{q} = [a_1, \ldots ,a_i,\frac{r}{d}].  \\  [1.8mm] 

\end{array}$

\vspace{.1in}
\noindent \, $6.$ \ \ \ If \ \ $a_i= b_i + c_i$ \ and \ $\frac{s}{u} = [c_i, a_{i+1},\ldots
,a_n],$ \ then \ $\frac{p}{q} = [a_1, \ldots, a_{i-1}, b_i + \frac{s}{u}]$ 

 \ and \ \ $\frac{p}{q}  =  [a_1,\ldots,  a_{i-1}, b_i + c_i, a_{i+1}, \ldots, a_n] = [a_1,\ldots, 
a_{i-1}, b_i, 0, c_i, a_{i+1}, \ldots, a_n].$              

 } \end{lem}

\begin{rem}{\rm \  If a continued fraction  $[a_1, a_2, \ldots, a_n]$ has even length, then we can 
bring it to odd length via the last term transformations: 
 \[ [a_1, a_2, \ldots, a_n] = [a_1, a_2, \ldots, a_n -1, +1] \mbox{ \rm for } a_n >0  \mbox{ \rm
and } \]
 \[ [a_1, a_2, \ldots, a_n] = [a_1, a_2, \ldots, a_n +1, -1] \mbox{ \rm for } a_n <0. \ \ \ \ \ \
\]
 } \end{rem}

\noindent We shall  say that a continued fraction is {\it
termwise positive} ({\it negative})  if all the numerical terms in its expression are positive
(negative).

\begin{defn}{\rm \  A  continued fraction $[\beta_1, \beta_2, \ldots, \beta_m]$ is said to be in 
{\it canonical form} if it is termwise positive or negative and 
$m$ is odd.   }
\end{defn}

\noindent By Remark 5 above any termwise positive or negative continued fraction may be assumed to 
be in canonical form.  The main observation now is the following well-known fact about continued
fractions
  (the analogue of Proposition 2).  

\begin{prop}{ \  Every  continued fraction $[a_1, a_2, \ldots, a_n]$ can be transformed to a 
unique canonical form with sign generically equal to the sign of the first non-zero term.
 } \end{prop}  

\noindent {\em Proof.}  Let $\frac{p}{q} = [a_1, a_2, \ldots, a_n]$  and suppose that the
$a_j$'s are not all of the same sign. Let $a_{i-1}, a_i$ be the first pair of adjacent $a_j$'s of
opposite sign, with $a_{i-1} > 0$.  We point out that the signs of $a_{i+1}, \ldots,  a_n$ are
irrelevant. We will show that $\frac{p}{q} = [\beta_1, \beta_2, \ldots, \beta_m],$ where all
$\beta_i$'s are positive or all negative. 
We do the same arithmetic operations to the continued fraction $[a_1, a_2, \ldots, a_n],$ as for
rational tangles and we check the results. Indeed, we have: 

\vspace{.15in}

\noindent $\begin{array}{rcl} 
\frac{p}{q} & = & [a_1, a_2, \ldots, a_n]  \\ [1.8mm] 

 \ & = &  [a_1, \ldots, a_{i-2}, (a_{i-1}-1)+1, -1+(a_i+1), a_{i+1}, \ldots, a_n]  \\    [1.8mm] 

 \ & \stackrel{L.6(6)}{=} &  [a_1, \ldots, a_{i-2},  (a_{i-1}-1), 0, +1, -1+(a_i+1), a_{i+1},
\ldots,
                           a_n]  \\    [1.8mm] 

 \ &\stackrel{L.6(6)}{=} &  [a_1, \ldots, a_{i-2}, (a_{i-1}-1), 0, +1, -1 + \frac{r}{l}],   \\ 
\end{array}$

\vspace{.15in}

\noindent where \ $\frac{r}{l} = [a_i+1, a_{i+1}, \ldots, a_n].$ \ This is transformed to

\vspace{.15in}

\noindent $\begin{array}{rcl} 
\frac{p'}{q'} & = &   [a_1, \ldots, a_{i-2},
(a_{i-1}-1), +1, -(a_i+1), -a_{i+1}, \ldots, -a_n] \\ 

 \ & = & [a_1, \ldots, a_{i-2}, (a_{i-1}-1), +1, -\frac{r}{l}].  \\    [1.8mm] 
\end{array}$

\vspace{.15in}

\noindent In order to show now that \ $\frac{p}{q} = \frac{p'}{q'}$ \ it suffices to show the
arithmetic equality 

\[ [0, +1, -1 + \frac{r}{l}] = [+1, -\frac{r}{l}]  \  \Longleftrightarrow \ 
\frac{1}{+1 + \frac{1}{-1 + \frac{r}{l}}} = +1  -\frac{l}{r},\]

\noindent which is indeed valid.  There is a similar identity for $a_{i-1} < 0$.  Notice that the
sum of the absolute values of the entries of the continued fraction $\frac{p'}{q'}$ is reduced by
one. So, proceeding by induction, we eliminate in the continued fraction all entries with negative
sign. Notice also that the sign of
$a_{i-1}$ and thus of
$a_1,$ if $a_1 \neq 0,$ dominates the above calculations.  As in the case of rational tangles  
(Proposition 2) there is one exception to this rule, namely when the continued fraction begins 
with an alteration of $+1$ and $-1$. More precisely, if $\frac{P}{Q} = [+1, -1, \frac{p}{q}],$
then 
$\frac{P}{Q} =  \frac{q}{q - p},$ and the sign of $\frac{P}{Q}$ is opposite to the sign of
$\frac{p}{q}.$ If $\frac{P}{Q} = [+1, -1, +1, -1, \frac{p}{q}],$ then 
$\frac{P}{Q} =  \frac{p - q}{p},$ and the sign of $\frac{P}{Q}$ is same as the sign of
$\frac{p}{q},$ and if $\frac{P}{Q} = [+1, -1, +1, -1, +1, -1, \frac{p}{q}],$ then 
$\frac{P}{Q} =  \frac{p}{q}.$ There are analogous identities for alterations of $-1$ and  $+1$. 
 Finally, by Remark 5, the index $m$ of the last term of the continued fraction $[\beta_1, \beta_2,
\ldots, \beta_m]$ can be assumed to be odd, and the uniqueness of the
final continued fraction follows from  Euclid's algorithm. This completes the proof. 
 $\hfill \Box $ 
\bigbreak

Another interesting fact about continued fractions is that any positive continued fraction
can be written as a continued fraction with even integer denominators, see \cite{Sie}. 
Note that, by Lemma 6(4), this fact can be extended to negative continued fractions. 
Siebenmann  \cite{Sie} uses this observation for finding an obvious Seifert surface spanning a 
given rational knot. 
\bigbreak

\noindent {\bf Matrix interpretation for continued fractions.} We now give a
way of calculating continued fractions via $2\times 2$ matrices (compare with \cite{Fr}, \cite{Ko},
\cite{Sie}, \cite{BZ}).  Let $[a_1, a_2, \ldots, a_n] = p/q.$ We correspond $p/q$  to the vector
$\left(
\begin{array}{cc}
     p \\
     q
\end{array}
\right)$ and we let $M(a_i) =  \left(
\begin{array}{cc}
     a_i & 1 \\
     1 & 0
\end{array}
\right) $ and $v = \left(
\begin{array}{cc}
     1 \\
     0
\end{array}
\right).$ Then, in this notation we have: \[ [a_1, a_2, \ldots, a_n]
= M(a_1) M(a_2) \cdots M(a_n) \, v. \]

\bigbreak
\noindent  {\bf Infinite tangles.} Before closing this section we push the analogy to periodic 
infinite tangles and imaginary tangles. It is a classic result that to every real number $r$
corresponds a unique continued fraction $[a_1, a_2,\ldots]$ that converges to $r,$  such that the
$a_i\in \ZZ$ and $a_i>0$ for all $i>1$ (see for example \cite{Kh}).  It is easy to see that we
could have instead 
 the $a_i$'s either all positive or all negative. This continued fraction is finite if 
$r$ is rational and infinite if $r$ is irrational. 
\smallbreak

 Further, it was proved by Lagrange that an
irrational number is quadratic (i.e. it satisfies a quadratic equation with integer coefficients) 
if and only if it has a continued fraction expansion which is periodic from some point onward.
(See 
\cite{Kh}, \cite{O}.)  Let  \, $\alpha \chi^2 =  \beta \chi +\gamma,$ \, be a quadratic equation
with integer coefficients and $\alpha\neq 0.$ The solutions
$\chi, \chi'$ will be either both real or both complex conjugates. 
 If the roots are real irrationals 
  we can find the periodic continued fraction expansion of one of the two
(the greater one, say $\chi$) by solving the equation \, $\chi = a_1 + \frac{1}{{\chi}_2},$ \,
where the number 
 \, ${\chi}_2 = \frac{1}{\chi - a_1} >1$ \, is irrational. We continue solving a similar  equation
for $\chi_2,$ and so on, until we obtain 
$\chi = [a_1, a_2, \ldots, a_k, \overline{b_1, b_2, \ldots, b_n}],$ where the bar marks the period
of the continued fraction. For example, the golden ratio is the positive  root of the equation 
$\chi^2 =  \chi +1,$ which gives rise to the infinite continued fraction $[1,1,1, \ldots].$  
 For the continued fraction expansion of the root $\chi'$ we know the following
remarkable theorem of Galois (also implicit in the work of Lagrange,  see for example
\cite{O}):  If $\chi >1$ is a quadratic irrational number and we have that $-1 < \chi' < 0,$ then
the  continued fraction expansion of $\chi$ is purely periodic.  Let $\chi = [\overline{a_1, a_2,
\ldots, a_n}]$ for $a_1, a_2,\ldots, a_n$ positive integers  
 and let  $\psi = [\overline{a_n, a_{n-1}, \ldots,
a_1}]$ be the continued fraction for $\chi$ with the period reversed. Then $-\frac{1}{\psi} =
\chi'$ is the conjugate root of the quadratic equation satisfied by $\chi.$ 
\smallbreak

 It is interesting to look at the relations of the above continued fractions
and corresponding infinite tangles. According to the above, each non-rational real number
(algebraic or transcendental) can be associated to an infinite tangle  $[[a_1], [a_2], [a_3], \ldots],$
all the approximants of which are rational tangles. A quadratic irrational number $\chi$ will be
associated to an infinite periodic rational tangle. This demonstrates a fractal pattern.  If  the
tangle for
$\chi$ is purely periodic, i.e. a tangle of the form $[\chi] = [\overline{[a_{1}], [a_{2}], \cdots
[a_{n}]}]$ then its conjugate will correspond to the $90^0 - $ rotation of this tangle with the period
reversed. In Figure 15 we illustrate the tangle for the golden ratio: 

$$[\frac{1+ \sqrt{5}}{2}] = [1] + \frac{1}{[1] + \frac{1}{[1] +
\frac{1}{[1]+\cdots}}}.$$

$$ \picill7.3cmby1.4in(ratl15) $$

\begin{center}
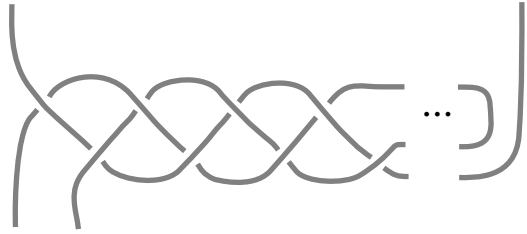
{Figure 15 - The tangle of the golden ratio $\frac{1+ \sqrt{5}}{2}$}
\end{center}
\vspace{3mm}

 Suppose now that  the quadratic equation  \, $\alpha \chi^2 =  \beta \chi +\gamma,$ \,
does not have real roots.  In this case we cannot apply the above algorithm for
obtaining an infinite continued fraction, whose limit value is well-defined. Yet we can write a
formal solution as an infinite continued fraction with rational entries, in the following way: 

$$ \chi^2 = \frac{\beta}{\alpha}\, \chi + \frac{\gamma}{\alpha} \ \ \ \Longrightarrow \ \ \ \chi =
\frac{\beta}{\alpha} +
\frac{\frac{\gamma}{\alpha}}{\chi}  \ \ \ \Longrightarrow \ \ \  \chi  = 
\frac{\beta}{\alpha} +
\frac{\frac{\gamma}{\alpha}}{\frac{\beta}{\alpha} + \frac{\frac{\gamma}{\alpha}}{\chi}} \, = \,  
\frac{\beta}{\alpha} +
\frac{1}{\frac{\beta}{\gamma} + \frac{1}{\chi}}.$$ 

\noindent Thus, with repeated iterations we obtain for $\chi$ the infinite purely periodic 
formal continued fraction with rational terms: 

$$ (\frac{\beta}{\alpha}) + \frac{1}{(\frac{\beta}{\gamma}) + \frac{1}{(\frac{\beta}{\alpha})
+ \frac{1}{(\frac{\beta}{\gamma}) +  \frac{1}{(\frac{\beta}{\alpha}) +  \cdots}}}}
\, = \, [\overline{\frac{\beta}{\alpha}, \frac{\beta}{\gamma}}].$$

\noindent The finitely iterated fraction values must oscillate in some set of values (possibly
infinite), and we  have behaviours of great  complexity related to the powers of the
complex number solutions.  In this form we can insert the rational tangles 
$[\beta/\alpha]$ and $[\beta/\gamma]$ into the places of horizontal and vertical twists
respectively  of the standard form of rational tangles illustrated in Figure 9 (where we have
previously restricted ourselves to integer and vertical tangles). 
The continued fraction form of the rational tangles 
$[\beta/\alpha]$ and $[\beta/\gamma]$ is found by writing out
the fractions $\beta/\alpha$ and $\beta/\gamma$ as continued fractions. 
 The  result is a sequence of
generalized continued fraction tangles that are not (even in the finite approximations) 
neccessarily rational. We shall call such tangles `imaginary'. 

\smallbreak
\noindent   For example, consider the equation 
\, $\chi^2 = \chi - 2.$ \, This has roots 
$\chi = \frac{(1+ \sqrt{-7})}{2}$ and $\chi' = \frac{(1- \sqrt{-7})}{2}.$ 
According to the above we can set up an infinite imaginary tangle with corresponding equation \,
$[\chi] = [\overline{[1] , \frac{1}{[-2]}}].$ \, We leave
it as an exercise for the reader to investigate $[\chi]$ and its finite
approximations. The finite approximations go chaotically through an
infinite set of fraction values.
Certainly $[\chi]$ deserves the name $[\frac{(1 + \sqrt{-7})}{2}].$
This is a case of using a rational insertion in the pattern of the
continued fraction forms. Another example is  $[\psi] = [\overline{[1], [-1]}]$ for
which the  corresponding formal infinite continued fraction is $[1, -1, 1, -1, \ldots].$ 
This leads to the equation \, $\psi = 1 + \frac{1}{-1 + \frac{1}{\psi}}$ \, 
 and to the quadratic equation \, $\psi^{2} = \psi - 1$ with roots $\psi = (1
\pm \sqrt{3}i)/2.$ \, 
The approximating fractions oscillate through the values $1, \ 1+ \frac{1}{-1} = 0$, \ 
$1 + \frac{1}{-1 + \frac{1}{1}} = \infty$ \, with period three. Notice that the periodic continued
fraction $[1, -1, 1, -1, \ldots]$ does not satisfy the conditions for convergence to a real
number.  Finally, another interesting example is the tangle \, $[i] = [\sqrt{-1}].$ \,
Here $i$ is a root of the quadratic equation $\chi^{2}+1 =
0,$ \, so \, $\chi  = -
\frac{1}{\chi}.$  Thus, the elemental imaginary  tangle  satisfies the equation \, $[\sqrt{-1}] =
-\frac{1}{[\sqrt{-1}]}.$  Since $-\frac{1}{T}$ is represented by the rotation $T^r,$ we see that  \,
$[\sqrt{-1}] = [\sqrt{-1}]^r.$ \, This is illustrated by the infinite tangle in Figure 16.

$$ \picill8cmby2.7in(ratl16) $$

\begin{center}
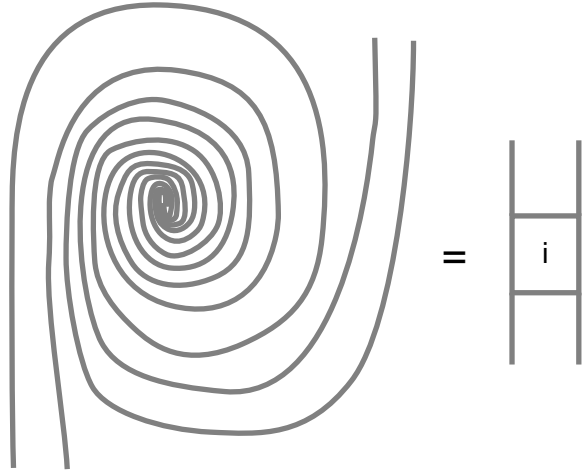
{Figure 16 - The tangle of the square root of -1}
\end{center}
\vspace{3mm}

\section{The Proof of the Classification Theorem}

 Let $T$ be a rational tangle in twist form:

\[ T = [s_k] +( \cdots + (\frac{1}{[r_3]} * ([s_1] + (\frac{1}{[r_1]} * [s_0] * \frac{1}{[r_2]}) +
[s_2]) * \frac{1}{[r_4]}) + \cdots )+ [s_{k+1}].\] 

\begin{defn}{\rm \ We define {\it the fraction of $T, \ F(T)$,}  to be the rational number  

\[ F(T) = s_k +( \cdots + (\frac{1}{r_3} * (s_1 + (\frac{1}{r_1} * s_0 * \frac{1}{r_2})
+ s_2) * \frac{1}{r_4}) + \cdots )+ s_{k+1},\] 

\noindent if $T \neq [\infty],$ and $F([\infty]) := \infty = \frac{1}{0},$ as a formal expression, 
where the arithmetic operation `$*$' is  defined via 
\[  x * y := \frac{1}{\frac{1}{x} + \frac{1}{y}}.  \]
 } \end{defn}

\noindent  For example we have: \ $F([0])=0, \ F([\pm1])=\pm1, \  F([\pm k])=\pm k, \
F(\frac{1}{[\pm k]}) =\frac{1}{\pm k}.$ Also, \ 
$F(([3] + (\frac{1}{[5]} * [6] * \frac{1}{[2]}) + [-4]) =  3 + \frac{1}{5 + \frac{1}{6} + 2} +
(-4).$ 

 \begin{lem}{ \ Let $T$ be a rational tangle in twist form and $C$ its continued fraction
form. Then $F(T) = F(C).$
 } \end{lem}
\noindent {\em Proof.} We observe first that, by Definition 8, the operation `$*$' is commutative. 
Also it is associative, since   $(a*b)*c = a*(b*c) = \frac{1}{\frac{1}{a} + \frac{1}{b} +
\frac{1}{c}}.$ 
 Thus, for the operations
`$+$' and `$*$' we have the identities: $F([n] + T ) = F(T + [n])$ \ and \  
$F(\frac{1}{[n]} * T) = F(T * \frac{1}{[n]}).$   For $T$ now with an expression as
above we have
\[ F(T) = s_k + \cdots  + \frac{1}{(r_3 + \frac{1}{(s_1 + \frac{1}{ (r_1 + \frac{1}{s_0} + r_2)} +
s_2)} + r_4)} + \cdots + s_{k+1}.  \]
\noindent  On the other hand we have from Lemma 3 that 
\[ C = ( \cdots ((( [s_0] * \frac{1}{[r_1 + r_2]}) +
[s_1 + s_2]) * \frac{1}{[r_3 + r_4]}) + \cdots )+ [s_k + s_{k+1}]). \] 
\noindent Thus 
\[ F(C) = (s_k + s_{k+1}) +  \cdots  + \frac{1}{(r_3  + r_4) + \frac{1}{(s_1 + s_2) + \frac{1}{ (r_1
+ r_2) + \frac{1}{s_0}}}}  = F(T).  \]
 $\hfill \Box$

\begin{rem}{\rm \  It follows from the above that:    
\[ \mbox{If} \ T= [a_1]+ \frac{1}{[a_2]+\cdots + \frac{1}{[a_{n-1}]
+\frac{1}{[a_n]}}} \ \ \mbox{then} \ \ F(T) = a_1+ \frac{1}{a_2+\cdots + \frac{1}{a_{n-1}
+\frac{1}{a_n}}},  \] and  this can be taken as the definition of $F(T).$ 
 } \end{rem}

 \begin{lem}{ \ Let $T= [[a_1], [a_2], \ldots, [a_n]]$ be a rational tangle in continued fraction
form. Then the tangle fraction has the following properties.

\vspace{.15in}

\noindent $\begin{array}{lrcl} 
 
 1. & F(T + [\pm 1]) &  = & F(T) \pm 1,   \mbox{ \ and \ }  F(T \pm [k]) = F(T) \pm k, \\   
[1.8mm] 

 2. & F(\frac{1}{T}) &  = & \frac{1}{F(T)},   \\  [1.8mm] 

3. & F(-T) &  =  & -F(T), \\  [1.8mm] 

4. & F(T * [\pm 1]) &  = & F(T) * (\pm 1),  \mbox{ \ and \ }  F(T* \frac{1}{[\pm k]}) = \frac{1}{\pm
k  + \frac{1}{F(T)}}, \\    [1.8mm] 

5. & \mbox{If \ } R &  =  & [[a_{i+1}],\ldots ,[a_n]], \mbox{ \ then \ } F(T)=[a_1, \ldots ,a_i,
F(R)],  \\  [1.8mm] 

6. & \mbox{If \ } a_i &  =  &  b_i + c_i  \mbox{ \ and \ }  S = [[c_i], [a_{i+1}],\ldots, [a_n]], 
\\  [1.8mm] 

 & \mbox{then \ } F(T) &  =  &   [a_1, \ldots, a_{i-1}, b_i + F(S)] =   
[a_1,\ldots, a_{i-1}, b_i, 0, F(S)]. \\  
\end{array}$
 } \end{lem}

\noindent {\em Proof.} Immediate from Lemmas 4, 5 and 6.  $\hfill \Box $

\bigbreak
 \noindent It follows from Lemma 8(2) that \ $F( \frac{1}{\frac{1}{T}} )=
F((T^{r})^{r})=F(T).$ 

\begin{lem}{ \ If  $T$ rational, then $F(T^{hflip}) =  F(T) = F(T^{vflip}).$
 } \end{lem}

\noindent {\em Proof.} We prove the first equality; the proof of the second one is completely
analogous. As for Lemma 2, we proceed by induction. The statement is true for the tangles
$[0], \ [\infty], \ [\pm 1],$ and assume it is also true for any rational tangle $R$ with $n$
crossings, i.e. $F(R) = F(R^{hflip}).$  By Remark 1, we only need to show that the statement is
valid for the  tangles $F~=~R~+~[\pm 1]$ and $L~=~R~*~[\pm 1].$ Indeed, for  $F^{hflip}$ and
$L^{hflip}$ we have:
\smallbreak
\noindent $F(F^{hflip}) = F((R+[\pm 1])^{hflip}) = F(R^{hflip}+[\pm 1]) \stackrel{L.8(1)}{=}
F(R^{hflip})\pm 1 \stackrel{induction}{=}$ 
\smallbreak
\noindent 
$F(R)\pm 1 \stackrel{L.8(1)}{=} F(R+[\pm 1]) = F(F), $ \ and
\smallbreak
\noindent $F(L^{hflip}) = F(R*[\pm 1])^{hflip} = F([\pm 1] * R^{hflip})
\stackrel{L.8}{=}  \pm 1 * F(R^{hflip}) \stackrel{induction}{=}$ 
\smallbreak
\noindent
$\pm 1 * F(R) \stackrel{L.7}{=}  F(R) * \pm 1 \stackrel{L.8(4)}{=} F(R*[\pm 1]) = F(L). $ 
 $\hfill \Box $

 \begin{lem}{ \ Let $T$ be a rational tangle in continued fraction
form and $T'$ its canonical form. Then $F(T) = F(T').$
 } \end{lem}

\noindent {\em Proof.} Direct consequence of the proofs of Propositions 2 and 3.  $\hfill \Box $

\bigbreak

We will show next that two alternating rational tangles are isotopic if and only if
they differ  by a finite sequence of flypes. Diagrams for knots and links are represented on the
surface of a two dimensional sphere and  then notationally on a plane for  purposes of illustration.
A {\it pancake flip} of a diagram is a  diagram obtained by picking up the diagram, turning it by
$180^0$ in  space and then replacing it on the plane. 

Abstractly we know that a diagram  and its
pancake flip are isotopic by Reidemeister moves. In fact, as we illustrate in Figure 17, a pancake
flip is a composition of $S^2$-isotopies, planar isotopy and a flype.  (By an {\it $S^2$-isotopy}
we mean the sliding of an arc around the back of the sphere.) 
To see this, note first that we can assume without loss of
generality that  we can isolate one crossing at the
`outer edge' of the diagram in the plane and decompose the diagram into this crossing and a
complementary tangle. I.e. the diagram in question is of the form 
$N([\pm 1] + R)$ for some tangle $R$ not necessarily rational.  In order to place the
diagram in this form we only need to use isotopies of the diagram in the plane. 
Thus, a pancake flip
is a composition of flypes up to $S^2$-isotopies, but it is convenient to have this move on
diagrams articulated directly.

$$\vbox{\picill5.5inby2.75in(ratl17)  }$$

\begin{center}
{ Figure 17 - Pancake flip} 
\end{center}
\vspace{3mm}

\begin{prop}{ \ Two alternating rational tangles on $S^2$ are isotopic if and only if they differ
by a finite sequence of rational flypes. 
 } \end{prop}

\noindent {\em Proof.} Let $T$ be a $2$-tangle contained in a $3$-ball in $S^3.$ By
shrinking the complementary $3$-ball to a point we may view it as a rigid vertex attached to the
tangle, see Figure 18. Thus, the {\it vertex closure} $V(T)$ is associated to the tangle $T$ in a
natural way. Note that $V(T)$ is an amalgamation of the numerator closure and the denominator
closure of $T,$ as defined in the introduction. An isotopy of $2$-tangles fixes their endpoints, so 
it can be considered as an isotopy of their vertex closures.

$$ \picill5.1cmby1.55in(ratl18) $$

\begin{center}
{Figure 18 - Vertex closure}
\end{center}
\vspace{3mm}
  
\noindent In \cite{SuT}, end of Section 1 it is argued that the solution to the Tait conjecture for
alternating knots implies that the flyping conjecture is also true for vertex closures of
alternating $2$-tangles and thus true for alternating $2$-tangles, see also \cite{Sa}. We shall
assume the Tait flyping conjecture for vertex closures of alternating  rational tangles and we
shall derive from this the flyping conjecture for alternating rational tangles. 

\smallbreak
 Let $T$ be an alternating rational tangle diagram. We consider all possible flypes on $V(T).$ 
 If a flype does not involve the rigid vertex of the closure then it is a tangle flype, thus by
Corollary 1 a rational flype, and so there is nothing to show. Consider now a flype that contains
the rigid vertex. We will show that such a flype can be reconfigured as the composition of a
pancake flip with a flype of a subtangle of the tangle
$T.$ Thus, up to a pancake flip, all flypes can take place on the tangle without involving the
vertex. 

Indeed, the region of a flype can be enclosed by a simple closed curve on the plane, that 
intersects the tangle in four points. Hence, a flype that involves the rigid vertex can only fall
into one of the two cases for $T:$ either $T=P+[\pm1]+R$ or $T=P*[\pm1]*R.$ 
Figure 19 illustrates for the first
case how to avoid to flype the rigid vertex up to a pancake flip.  Note
that we have shaded one arc of the rigid vertex darker, in order to make the isotopies easier to follow. 
The second case for $T$ follows from the first one by a $90^0 - $rotation on the plane.

$$ \picill8inby3.3in(ratl19) $$

\begin{center}
{Figure 19 - Vertex flype analysis} 
\end{center}
\vspace{3mm}

  Let now  $T$ and $S$ be two isotopic alternating rational tangles and let 
$V(T)$ and $V(S)$ be their vertex closures. By \cite{SuT} we have that $V(T)$ and $V(S)$ are 
related by a sequence of flypes. From the above reasoning it can be assumed that, up to a pancake
flip, these flypes all leave the rigid vertex fixed, hence they are tangle flypes. Now, the
horizontal pancake flip induces a horizontal flip and the vertical pancake flip induces a vertical
flip on the rational tangle. These, by Lemma 2, are isotopic to the original rational tangle.
Thus, all steps above are tangle isotopies. Finally, by Corollary 1, 
tangle flypes on rational tangles have to be rational. This completes the proof.  
$\hfill \Box $

\begin{cor}{ \ It follows from Lemma 1 and Proposition 4 that two isotopic rational tangles with 
all crossings of the same type will be twist forms of the same canonical form.}
\end{cor}   

 \begin{lem}{ \ Two rational tangles that differ by a rational flype have the same fraction. 
 } \end{lem}
\noindent {\em Proof.} Let $T$ and $S$ be two rational tangles that differ by a flype with respect
to a rational subtangle $t.$ The flype will have one of the algebraic expressions:  $[\pm 1] +  t
\sim t^{hflip} +[\pm 1]$ \, or \, $[\pm 1] * t \sim t^{vflip} * [\pm 1].$ By Lemma 9 \, 
$F(t^{hflip}) = F(t)$ \, and  \, $F(t^{vflip}) = F(t),$ \, and by Lemma 7 \, 
 $F([\pm 1] + t) = F(t +[\pm 1])$ \, and \, $F([\pm 1] * t) = F(t * [\pm 1]).$ Finally, 
by Corollary 1 \, $t$ is a rational truncation of $T,$ and Lemmas 5 and 6 tell us that continued
fractions of rational tangles and arithmetic continued fractions agree on truncations. Thus,  
we obtain $F(T)=F(S).$ 
 $\hfill \Box $

\begin{th}{ \ The fraction is an isotopy invariant of rational tangles.
 } \end{th}
\noindent {\em Proof.} Let $T, S$ be two isotopic rational tangles in twist form. By Lemma 3 and
Proposition 1 the tangles $T, S$ can be isotoped to two rational tangles $T', S'$ in continued
fraction form, and by Lemma 7 we have $F(T) = F(T')$ and $F(S) = F(S').$ Further, by Proposition
2 the tangles $T', S'$  can  be isotoped to two alternating rational tangles $T'', S''$ in
canonical form, and by Lemma~10  we have
$F(T') = F(T'')$ and $F(S') = F(S'').$ Finally, by Proposition 4 the tangles $T'', S''$ will differ
only by rational flypes, and by Lemma 11 we have $F(T'') = F(S'').$ Thus $F(T) = F(S),$ and this
ends the proof of the theorem.
 $\hfill \Box $

\begin{th}{ \ Two rational tangles with the same fraction are isotopic.
 } \end{th}
\noindent {\em Proof.}  Indeed, let $T = [[a_1], [a_2], \ldots, [a_n]]$ and
$S = [[b_1], [b_2], \ldots, [b_m]]$ be two  rational tangles with 
$F(T) = F(S) = \frac{p}{q}.$  We bring $T, S$ to their canonical forms
$T' = [[\alpha_1], [\alpha_2], \ldots, [\alpha_k]]$ and $S' = [[\beta_1], [\beta_2], \ldots,
[\beta_l]]$  respectively. From Theorem 2 we have $F(T') = F(T) = F(S) = F(S')
= \frac{p}{q}.$  By Proposition 3, the fraction $\frac{p}{q}$ has a unique continued fraction 
expansion in canonical form, say
$\frac{p}{q} = [\gamma_1, \gamma_2, \ldots, \gamma_r].$  This gives rise to the  alternating
rational tangle in canonical form
$Q = [[\gamma_1], [\gamma_2], \ldots, [\gamma_r]],$ which is uniquely determined from the vector of
integers  $(\gamma_1, \gamma_2, \ldots, \gamma_r).$ We claim that  $Q=T'$ (and similarly $Q=S'$).  
Indeed, if this were not the case we would have the two different continued fractions  in canonical
form giving rise to the same rational number: $[\alpha_1, \alpha_2, \ldots, \alpha_k] =
\frac{p}{q} = [\gamma_1, \gamma_2, \ldots, \gamma_r].$ But this contradicts the uniqueness of the
canonical form of continued fractions  (Proposition 3).
 $\hfill \Box $

\bigbreak 
\noindent {\bf Proof of Theorem 1.} \ Theorems 2 and 3 show that two rational tangles are isotopic
if and only if they have the same fraction, yielding the proof of Theorem 1 as a corollary. 
 \hfill Q.E.D.

\bigbreak
We conclude this section with some comments.

\bigbreak

\noindent {\bf Note 2} \  It follows from Theorem 1 
that if $T = [[a_1],[a_2],\ldots ,[a_n]]$ is a rational tangle in continued fraction form, and if 
$\frac{p}{q} = [a_1,a_2,\ldots,a_n]$ is the evaluation of the corresponding arithmetic continued
fraction then, without ambiguity, we can write $T = [\frac{p}{q}].$ Thus, rational numbers are
represented bijectively by rational tangles, their negatives are represented by the mirror images 
and their inverses by the inverses of the rational tangles. 

Moreover, adding integers to a
rational number corresponds to adding integer twists to a rational tangle, but sums of non-integer
rational numbers do not correspond to the rational tangles of the sums. Such sums go beyond the
rational tangle category; they give rise to `algebraic tangles'.   We call a tangle  {\em
algebraic} if it can be obtained by substituting rational tangles into an algebraic expression
generated from some finite set of variables by tangle addition and inversion.

Further, given a rational tangle
in twist or  standard form, in order to bring it to its canonical  form one simply has to
calculate its fraction and express it in canonical form. This last one gives rise to an
alternating tangle in canonical form which, by Theorem 1, is isotopic to the initial one. For
example, let $T = [[2],[-3],[5]].$ Then $F(T) = [2, -3, 5] =
\frac{23}{14}.$ But $ \frac{23}{14} = 
 [1,1,1,1,4],$ thus $T \sim [[1],[1],[1],[1],[4]],$ and this last tangle is the canonical form of
$T.$ 

\bigbreak
From the uniqueness of the canonical form of a continued fraction we  also have that: 

\begin{cor}{ \  The canonical form of a rational tangle is unique.  }
\end{cor}

\begin{cor}{ \ Rational tangles in canonical form have minimal number of crossings. }
\end{cor}   
\noindent {\em Proof.} \  Let $T''$ be a rational tangle in  canonical form and let ${\cal T}$ be the
set of all rational tangles in twist form with canonical form the tangle $T''.$ By Corollary 4, for 
each element of ${\cal T}$ the canonical form $T''$ is unique. 
Let now $T\in {\cal T}$ be a  rational tangle with $k$ crossings in twist form. By a sequence of 
flypes we bring $T$ to standard form $T' \sim T,$ and since flypes do not change the number of
crossings it follows that $T'$ has $k$ crossings. Note that $T'\in {\cal T}.$ We bring  $T'$ to its
canonical form 
$T'',$ and by the proof of Proposition 2, $T''$ will have less crossings than $T'.$ 
$\hfill \Box $

\begin{cor}{ \ Alternating rational tangles have minimal number of crossings. }
\end{cor}   
\noindent {\em Proof.} \  Indeed, if an alternating rational tangle is in twist form then by a
sequence  of flypes we bring it to canonical form, which by Corollary 5 has  minimal number of
crossings. And since flypes do not change the number of crossings the assertion is proved.
$\hfill \Box $

\bigbreak

\section{The Fraction through Integral Coloring} 

In this section we show how to compute the fraction of a rational tangle by coloring the arcs 
of the tangle with integers. This section is self-contained and does not depend upon the
development of the fraction that we have already made. So we eliminate the need for using the
Tait conjecture in our proof of classification of rational tangles. 

We used the Tait conjecture to
show that if two alternating rational tangles are isotopic then their fractions are equal. Without
the Tait conjecture we showed that if they have same fraction they are isotopic. Here we get the
isotopy invariance by the definition of the fraction. Thus, in combination with the Sections 2 and
3 and Theorem 3 this section provides another elementary proof of the classification of rational
tangles. The coloring method explained here is special to rational tangles and some of their
generalizations. The coloring gives an efficient and reliable method for computing the fraction of a
rational tangle (and from this its canonical form). Along with producing the fraction, the coloring
itself is of interest and it can be used to investigate related colorings of the closures of the
tangle. (See for example \cite{DJP, P, P2}.) 
\bigbreak

\noindent  We shall use colors from either $\ZZ$ or from $\ZZ_n$ for some n. The
coloring rule is that if two undercrossing arcs colored
$\alpha$ and
$\gamma$  meet at an overcrossing arc colored $\beta$, then \, $\alpha + \gamma = 2\beta.$
See Figure 20.  We often think of one of the undercrossing arc colors as determined
by  the other one and the color of the overcrossing arc. Then one writes $\gamma = 2\beta -\alpha.$
It is easy to verify that this coloring method  is invariant under 
the Reidemeister moves in the following sense: Given a choice of coloring for  the tangle (knot),
there is a way to re-color it each time a Reidemeister move (or a flype) is performed, so that no
change occurs to the colors on the external strands of the tangle (so that we still have a valid
coloring).  This means that a  coloring potentially contains  topological information about a
knot or a tangle. 
\smallbreak
\noindent In coloring a knot (and also many non-rational tangles) it is usually necessary to reduce 
the colors to the set of integers modulo
$N$ for some modulus $N$. In Figure 20 it is clear that the color set ${\ZZ}/{3\ZZ} =
\{ 0,1,2 \}$ is forced for coloring a trefoil knot.

$$ \picill5.2inby1.5in(ratl20) $$

\begin{center}
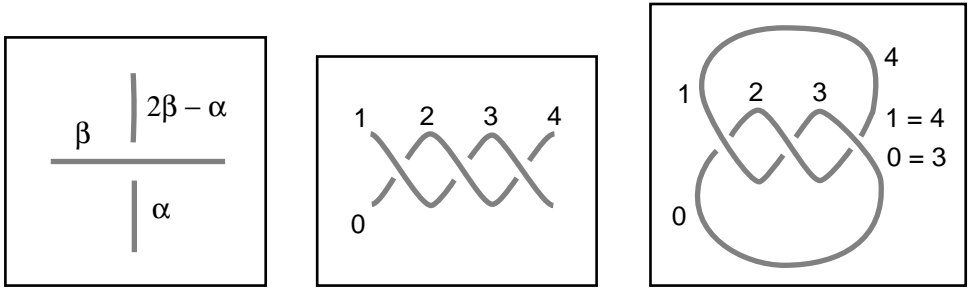
{Figure 20 - The coloring rule, integral and modular coloring } 
\end{center}

\vspace{3mm}

\noindent   When there exists a coloring of a tangle by 
integers, so that it is not necessary to reduce the colors over some modulus we shall say that the
tangle is {\it integrally colorable}.  
It turns out that {\it every rational tangle is integrally colorable:} Choose two colors for the
initial strands (e.g. the colors $0$ and $1$) and color the rational tangle as you create  it by
successive twisting.  We call the colors on the initial strands the {\it starting colors}. It is
important that we start coloring from the initial strands, because then the coloring propagates
automatically and uniquely. If one starts from somewhere else, one might get into an edge with an
undetermined color.

$$ \picill5.1inby2.8in(ratl21) $$

\begin{center}
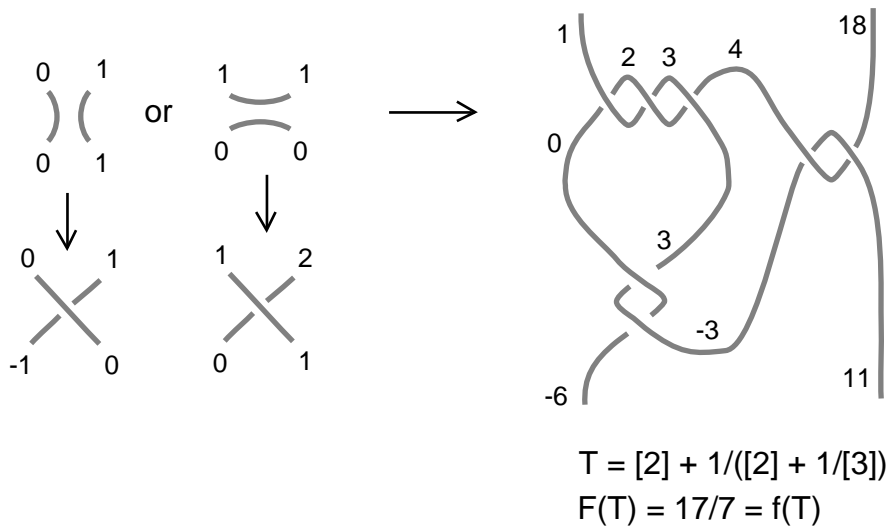
{Figure 21 - The starting colors, coloring rational tangles } 
\end{center}

\vspace{3mm}

\noindent The resulting colored
tangle now has colors assigned to its external strands at the northwest, northeast, southwest and
southeast positions. Let $NW(T)$,
$NE(T)$,
$SW(T)$ and
$SE(T)$ denote these respective colors of the colored tangle
$T$ and define the {\em color matrix of $T$}, $M(T)$, by the equation

$$M(T) =  \left[
\begin{array}{cc}
     NW(T)& NE(T) \\
     SW(T) & SE(T)
\end{array}
\right]. $$

\noindent We wish to extract topological information about the rational tangle $T$ from this 
matrix. Letting 

$$M =  \left[
\begin{array}{cc}
     a& b \\
     c & d
\end{array}
\right] $$

\noindent be a given color matrix we see at once from the above description of the coloring of a 
rational tangle that

$$M' =  \left[
\begin{array}{cc}
     na+k & nb+k \\
     nc+k & nd+k
\end{array}
\right]$$

\noindent will also be a color matrix for the given tangle. To see this replace each color
$\alpha$ by the color $n \alpha + k$ and note that if $\gamma = 2 \beta - \alpha$ then
$n \gamma + k = 2 n \beta + k - (n \alpha + k).$ Hence the new coloring is indeed
a coloring and the endpoints are replaced as indicated.   As a result of
this observation, we see that it is possible to set the starting colors  equal to $0$ and $1$ and
that this will change the color matrix by a sequence of transformations of the type $M
\longmapsto M'$ shown above.

\begin{th}{ \ 
 Let 

$$M =  \left[
\begin{array}{cc}
     a& b \\
     c & d
\end{array}
\right] $$

\noindent be a color matrix for an integrally colored tangle $T$.  Then  

 \begin{itemize}
\item[1.] \  $M$ satisfies the `diagonal sum rule': $a+d=b+c$.
\vspace{-.1in}
\item[2.] \  If $T$ is rational, then the quantity

$$f(T) := \frac{b-a}{b-d}$$ 

\noindent is a topological invariant associated with the tangle $T.$ 

\vspace{-.1in} 
\item[3.] \ $f(T+S) = f(T) + f(S),$ \ when there is given an integral coloring of a tangle $T+S$.
The colorings of $T$ and
$S$ are the restrictions of the coloring of $T+S$ to these subtangles. 
 
\vspace{-.1in}
\item[4.] \ $f(-\frac{1}{T}) = -\frac{1}{f(T)}$

\noindent  for any integrally colored $2$-tangle $T$ satisfying the diagonal sum rule.
\vspace{-.1in}
\item[5.] \  $f(-T) = -f(T)$ \ for any rational tangle $T$. Hence, 
\vspace{-.1in}
\item[6.] \ $f(\frac{1}{T}) = \frac{1}{f(T)}$  \ for any rational tangle $T$.
\vspace{-.1in}
\item[7.] \ $f(T) = F(T)$  \ for any rational tangle $T$.
\end{itemize} 

\noindent Thus the coloring fraction is identical to the arithmetical fraction defined earlier. 
 } \end{th}

\noindent We note that if $T$ is colored but not 
rational, we let $f(T)$ be defined by the same formula, but note that it may depend on the choice
of coloring. 
\bigbreak
\noindent {\em Proof.} It is easy to see that there are colorings  for $[0]$ and $[1]$ (see
Figure 21) so that
$f([0]) =
\frac{0}{1}$, 
$ f([\infty])=
\frac{1}{0}$, $f([1])=1.$ Hence property 7 follows by 3, 5 and induction. To see that the  diagonal sum
rule is satisfied for colorings of rational
 tangles, note that $a+d =
b+c$ implies that $d-c = b-a$ and $d-b = c-a.$  Then we proceed by induction  on the number of crossings
in the tangle. The diagonal sum rule is satisfied for colorings of the $[0]$ or $[\infty]$ tangle, since
the matrix for a coloring of such a tangle consists in two equal rows or two equal columns. Now  assume
that 

$$M =  \left[
\begin{array}{cc}
     a& b \\
     c & d
\end{array}
\right] $$

\noindent is a matrix for a coloring of a given tangle $T$ satisfying the diagonal sum rule. Then
it is easy to see that  $T + [1]$ has color matrix   

$$\left[
\begin{array}{cc}
     a& 2b-d \\
     c & b
\end{array}
\right] $$

\noindent and the identity $a+b = (2b-d) +c$ is equivalent to the identity $a+d=b+c$.  Thus the 
induced coloring on  $T+[1]$ satisfies the diagonal sum rule. The same argument applies to adding a
negative twist, as well as a twist on the left, bottom or top  of the tangle. Thus we
have proved by induction that the diagonal sum rule is satisfied for colorings of rational tangles. We
leave it as an exercise for the reader to prove the diagonal sum rule for any integrally colored
 $2$-tangle.  To show that $f(T) = (b-a)/(b-d)$ is a topological invariant of the tangle $T$
note that,  by  definition, the quantity  $f(T)$ is unchanged by the matrix transformations $M
\longmapsto M'$ discussed prior to the statement of this proposition. Thus, $f(T)$ does not depend upon
the choice of coloring for the rational tangle. Since, for any given coloring, $f(T)$ is a topological
invariant of the tangle with respect to that coloring, it follows that $f(T)$ is a topological
invariant of the tangle, independent of the choice of coloring used to compute it. 
\bigbreak

\noindent For proving property 3, suppose that $T$ has color matrix $M(T)$ and $S$ has color
matrix $M(S)$. Then for these to be the  restrictions from a coloring of $T+S$ it must be that the
right column of
$M(T)$ is identical with the left column of $M(S)$.  Thus

$$M(T) =  \left[
\begin{array}{cc}
     a& b \\
     c & d
\end{array}
\right], \ \ \ M(S) =  \left[
\begin{array}{cc}
     b& e \\
     d & f
\end{array}
\right], \ \ \ M(T+S) =  \left[
\begin{array}{cc}
     a& e \\
     c & f
\end{array}
\right]. $$

\noindent  Note that by the diagonal sum rule for $S$,  $b-d = e-f.$ Then 

$$f(T) + f(S) = \frac{b-a}{b-d} + \frac{e-b}{e-f} = \frac{b-a}{e-f} + \frac{e-b}{e-f}=
\frac{e-a}{e-f} = f(T+S).$$

\noindent This shows that $f(T)$ is additive with respect to tangle addition. Given $M(T)$ as
above, we have $M(-\frac{1}{T}) = M(T^r)$ given by the formula below: 

$$M(-\frac{1}{T}) =  \left[
\begin{array}{cc}
     b& d \\
     a & c
\end{array}
\right]. $$

\noindent Thus $$f(-\frac{1}{T}) = \frac{d-b}{d-c} = \frac{d-b}{b-a} = {-1}/{(\frac{b-a}{b-d})} =
-\frac{1}{f(T)},$$ 

\noindent and so property 4 is proved. The tangle $-T$ is obtained from the tangle $T$ by switching
all the crossings in $T$.  Let $T'$ be the tangle obtained from $T$ by reflecting it in a plane $P$
perpendicular to the plane on which the diagram of $T$ is drawn, as illustrated in Figure 22, i.e.
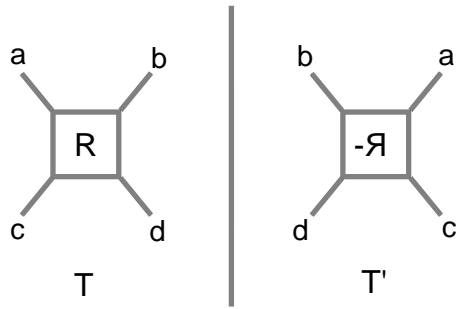
$T' := (-T)^{vflip}.$ We shall call $T'$ {\it the vertical reflect} of $T.$

$$ \picill2.75inby1.7in(ratl22) $$

\begin{center}
{Figure 22 - The vertical reflect of $T$}
\end{center}
\vspace{3mm}

\noindent It is then easy to see that a coloring of $T$ always induces a
coloring of $T'$ (the same colors that appear in $T$ will also appear in $T'.$) In fact,
if 
$$M(T) =  \left[
\begin{array}{cc}
     a & b \\
     c & d
\end{array}
\right]  \ \mbox{is a color matrix for} \ T, \ \mbox{then} \  M(T') =  \left[
\begin{array}{cc}
     b & a \\
     d & c
\end{array}
\right] $$
\noindent is the matrix for the induced coloring of $T'$. Therefore, using $a+d=b+c$, we have

$$f(T') = \frac{a-b}{a-c} = \frac{a-b}{b-d} = -\frac{b-a}{b-d} = -f(T).$$

\noindent By Lemma 2, $T'$ is isotopic to  $-T$ for rational tangles. So,
property 5 is proved. Property 6 follows from 4 and 5.  This completes the proof.
$\hfill
\Box$

\begin{rem}{\rm \  Rational tangles are integrally colorable, and it is easy to see that sums of
rational tangles are also integrally colorable. Also, it is easy to see that algebraic tangles are
integrally colorable (recall definition in Note 2).  At this writing, it is an open problem to
characterize integrally colorable tangles. The presence of a local knot, can keep a tangle from
being integrally colorable (by forcing the coloring into a specific modulus), but knotted arcs can
occur in integrally colorable tangles. For example, the non-rational algebraic tangle
$\frac{1}{[3]} +
\frac{1}{[2]}$ is integrally colorable and has a knotted arc in the form of the trefoil knot
(linked with another arc in the tangle). 
 } \end{rem}

\begin{rem}{\rm \ Note that if we have a tangle $T$ with color matrix
$$M(T) =  \left[
\begin{array}{cc}
     a & b \\
     c & d
\end{array}
\right],$$
\noindent we can subtract the color $a$ from all colors in the tangle, obtaining a new coloring with
matrix
$$M'(T) =  \left[
\begin{array}{cc}
     0 & b-a \\
     c-a & d-a
\end{array}
\right].$$
\noindent By the diagonal sum rule this has the form
$$M'(T) =  \left[
\begin{array}{cc}
     0 & a' \\
     b' & a'+b'
\end{array}
\right].$$
\noindent In thinking about colorings of tangles, it is useful to understand that one can always shift
one of the peripheral colors to the value zero. } \end{rem}

\begin{rem}{\rm \  Let $T$ be an $(m,n)$-tangle that is colored integrally, and suppose that
$a_1, a_2, \ldots, a_m$ are the colors from left to right on the top $m$ strands of $T$, and that
$b_1, b_2, \ldots, b_n$ are the colors from left to right on the bottom $n$ strands of $T$. Show
that 

\[ \Sigma_{i=1}^{m}(-1)^{i+1}a_i = \Sigma_{j=1}^{n}(-1)^{j+1}b_i \]

\noindent This is a generalization of the diagonal sum rule (pointed out to us by W.B.R. Lickorish
\cite{L}.)
\smallbreak

Consider, now, the knot or link $K = N(T).$ In order for the coloring of $T$ to be a coloring of
$K$, we then need that $a \equiv b$ and that $c \equiv d$. Since $a-b = c-d$ (by the diagonal sum
rule), we can take the coloring of $K$ to have values in ${\ZZ}/{D\ZZ}$ where $D = a-b$. This is an
example of a coloring of a knot occurring in a modular number system. This is more generally the
case, and one can always attempt to color a knot in ${\ZZ}/{Det(K)\ZZ},$ where $Det(K) =
|<K>(\sqrt{i})|,$  {\it the determinant of the knot}, where $<K>$ denotes the Kauffman bracket
polynomial of the knot $K.$  There are many fascinating combinatorial/topological problems related to coloring of
knots and tangles. 
 } \end{rem}

\begin{rem}{\rm \ View Figure 21 and note that the rational tangle $T = [2] + 1/([2] + 1/[3])$
with fraction $17/7$ is colored by starting with colors $0$ and $1$ at the generating arcs of the 
tangle and that {\it all the colors are distinct from one another as integers}. Furthermore, if
one takes the numerator closure $K = N(T)$ and colors in $\ZZ/{17\ZZ}$, the colors remain distinct in
this  modulus. This is not an accident! This is part of a more general conjecture about  coloring
alternating knots. See \cite{HK}. Here we prove the conjecture for rational knots and links. The
general result is stated below after a few preliminary definitions.
\bigbreak

If a crossing in a link diagram is regarded as the tangle $[+1]$ or $[-1]$ then it can be replaced 
by the tangle $[0]$ or the tangle $[\infty]$, maintaining the same outward connections with the
rest of the diagram. Such a replacement is called a {\it smoothing} of the crossing. A connected
link diagram is said to have a {\it nugatory crosssing} if there is a crossing in the link diagram
such that one of the smoothings of the diagram yields a disconnected diagram with two non-empty
components. In other words, at a nugatory crossing the diagram falls
apart into two pieces when it is smoothed in one of the two possible ways. We say that a diagram 
is {\it reduced} if it is connected and has no nugatory crossings. One can see easily that any
rational tangle diagram with no simplifiying Reidemeister one moves is a reduced diagram.

\begin{th}{ \  Let $T$ be a reduced alternating rational tangle diagram in twist form. Let $C(T)$ 
be any coloring of
$T$ over the integers. Then all the colors appearing on the arcs of $T$ are mutually distinct.
Furthermore, let $K = N(T)$ be the numerator closure of $T$ and suppose that the determinant of the
link $K$ is a prime number $p.$ Then for any coloring of $K$ in $\ZZ/{p\ZZ}$, all the colors on the arcs
of 
$K$ are distinct in $\ZZ/{p\ZZ}.$ In other words, if $v(K)$ denotes the number of crossings in the
diagram
$K$, then there will be $v(K)$ distinct colors in any coloring of the diagram $K$ in $\ZZ/{p\ZZ}.$
} \end{th}

\noindent {\em Proof.} The key to this proof is the observation that when one colors a reduced 
rational tangle starting with the integers $0$ and $1$ at the generating arcs, then all the
colors on the other arcs in the tangle are mutually distinct and increase or decrease in absolute
value  so that the largest colors in absolute value are the ones on the outer arcs of the tangle.
We have illustrated this phenomena in Figure 21. Note that in this figure the colors literally
increase as one goes through the first horizontal twist out to colors
$3$ and $4.$ Then we enter a sequence that is descending to $-1$ and $-6.$ The point to note is 
that this second sequence is genuinely descending and hence the sequence of numbers starting from
$-1$ and
$-6$ is ascending to $3$ and $4.$ The remaining twist sequence ascends to $11$ and $18.$ We leave 
it as an exercise for the reader to show by induction that this distinctness with maximal value
at the periphery holds for any reduced alternating rational tangle in twist form. 

Having checked
this property for tangles with starting values of $0$ and $1$ we can now assert its truth for
all colorings of the rational tangle by integers. All such colorings are obtained from the given
one by multiplying all colors by a non-zero constant or by adding a constant to each label in the
coloring. Distinctness and maximality is preserved by these arithmetical operations. Now consider 
the numerator closure
$K = N(T).$ It is not hard to see  (and we leave the proof for the reader) that if we start with 
colors
$0$ and $1$ at the generating arcs of the tangle, and if the resulting coloring has color matrix
$$M(T) =  \left[
\begin{array}{cc}
     a & b \\
     c & d
\end{array}
\right],$$
\noindent then $Det(K) = \pm (b-a).$ By the above discussion we can assume that $b$ and $a$ are 
the largest colors in absolute value on the diagram of $T$. Hence when we color $K$ in the modulus
$M = |Det(K)|$ we find that all the colors on $K$ are distinct in $\ZZ/{M\ZZ}.$ This proves that 
the chosen coloring for $K$ has the distinctness property. Now suppose that $N$ is a prime number
$p.$ Then $\ZZ/{M\ZZ} = \ZZ/{p\ZZ}$ is a field and hence the operation of multiplication of
colors by an non-zero element of $\ZZ/{p\ZZ}$ is invertible. It follows that all colorings
constructed from the given coloring by addition of a constant or multiplication by a non-zero
constant share in the distinctness property. Since these constitute all the non-trivial colorings
of $K$ over $\ZZ/{p\ZZ}$, the proof is complete.
$\hfill \Box$
\bigbreak

Theorem $5$  constitutes a proof, for rational knots and links, of a conjecture of Kauffman and 
Harary \cite{HK}. The conjecture states that {\it if $K$ is a reduced, alternating link diagram, and $K$ has
prime determinant $p$ then every coloring of the diagram $K$ in $\ZZ/{p\ZZ}$ has $v(K)$ distinct colors,
where $v(K)$ denotes the number of crossings in the diagram $K.$}  The conjecture has been independently
verified for rational knots and links and for certain related families of links in \cite{LP}.
 } \end{rem}

\begin{rem}{\rm \  Finally, we note that there is the following mapping

$$J: Color \, Matrices \longrightarrow \CC $$

\noindent induced via $$J(M(T)) := J
\left (\left[
\begin{array}{cc}
     a & b \\
     c & d
\end{array}
\right]\right)
:= (b-a)+i(b-d),$$

\noindent where  $Color \, Matrices$ denotes the set of color matrices satisfying the diagonal sum
condition. If $M$ is a color matrix, let $M^r$ be the color matrix obtained by rotating $M$
counterclockwise by $90^0$. Thus

$$M = 
\left [
\begin{array}{cc}
     a & b \\
     c & d
\end{array}
\right], \quad 
M^r = 
\left [
\begin{array}{cc}
     b & c \\
     a & d
\end{array}
\right].$$

\noindent Note that if $M=M(T),$ then $M^r=M(T^r),$ the matrix of the rotate of the tangle $T$.
Then it is easy to see that 

$$J(M^r) = i\cdot J(M).$$ 

\noindent Usually multiplication by $i$ is interpreted as a $90^0 - $ rotation of vectors in the
complex plane. With the equation 

$$J(M(T^r)) =  J(M(T)^r) = i\, J(M(T))$$ 

\noindent we see a new interpretation of $i$ in terms of $90^0$ rotations of tangles or matrices.
 } \end{rem}

We would like to conclude this section by a brief description of the fraction of rational tangles
through {\it conductance.} Conductance is a quantity defined in electrical networks as the inverse
of resistance. In \cite{GK1}  the conductance is defined as
a weighted sum of maximal trees in a graph divided by a weighted sum of maximal trees in an
associated graph, that is obtained by identifying the input and output vertices of the original
graph. This definition allows negative values for 
conductance and it agrees with the classical one, implying that in the resistance one would have to
consider also the notion of an amplifier. 

Conductance satisfies the law of parallel and series
connection as well as the star-triangle relation for appropriate values. Given a knot diagram one
can associate a graph, so that the Reidemeister moves on the knot diagram correspond to parallel
and series connection of resistances (Kirkhoff laws) and the star-triangle changes in the graph.
By defining the conductance on the knot diagram as the conductance on the corresponding graph one
shows that the conductance is an isotopy invariant of knots. The conductance of a rational tangle
turns out to be the numerical fraction of the tangle and from the above it does not depend on its
isotopy class.

\section{Negative Unity, the Group $SL(2,\ZZ)$ and Square Dancing}

The main result of this last section is integral to an illustrative game for the Conway Theorem on 
rational tangles. In this game (called `Square Dancing' by Conway) four people hold two ropes,
allowing the display of various tangles. The `dancers' are allowed to perform two basic moves 
called {\em turn} and {\em add}. Adding corresponds to an interchange of two dancers that adds one
to the corresponding tangle. Turning is a rotation of all four dancers by ninety degrees,
accomplishing negative reciprocation of the tangle.  We will show in this section that all rational
tangles can be produced by these operations, so the players can illustrate the classification
theorem.  

\bigbreak 

\noindent It is an interesting fact that the operations of {\it rotation and $+[1]$} generate
all rational tangles from the starting tangle of $[0]$. In order to see this, we generate the
operation $-[1]$ (which is the same as $+[-1]$) by iteration of the other two. Indeed, we have:

\begin{lem}{ \ The following identity holds for all rational tangles $x$.   

$$x-[1] = \frac{-1}{\frac{-1}{\frac{-1}{x}+[1]}+[1]}.$$ 

 } \end{lem}

\noindent {\em Proof.} The thing is that this identity holds for real numbers, thus showing that 
all rational numbers are generated by {\it negative reciprocation}  and {\it
addition of 1}. Since we know that arithmetical identities about rational tangles correspond to
topological identities the above identity is also valid for rational tangles.  This is the
arithmetic proof.  $\hfill \Box $
\bigbreak 
\noindent  Note that this property is equivalent to saying that 

$$ (r \circ (+1))^3 (x) = x, $$ 

\noindent where $r$ stands for the rotation operation, $+1$ for adding $[1]$, and $\circ$ for
composition of functions. That the three-fold iteration of $r \circ (+1)$ gives  the identity on
any tangle $T$ is illustrated in Figure 23, where we see that after applying $r \circ (+1)$ three
times to $T$, one of the tangle arcs can be isotoped to that the whole tangle is just a turned
version of the original.

$$ \picill8inby2.4in(ratl23) $$

\begin{center}
{Figure 23 - $(r \circ (+1))^3 = id$}
\end{center}
\vspace{3mm}

\noindent We also note that the statement of Lemma 12 can be modified for any $2$-tangle. 
 Now it reads $(r \circ (+1))^3 (T) = T^{r^2}.$ Figure 23 illustrates the general proof. 

\bigbreak

\noindent  In the header to this section, we advertized the
group $SL(2,\ZZ).$ The point is that Lemma 12  shows that the arithmetic of rational tangles  is
just isomorphic to the arithmetic of integer  $2 \times 2$ matrices of determinant equal to $+1$
(that being the definition of $SL(2,\ZZ)$.) The key point is the well-known fact that $SL(2,\ZZ)$ is
generated by matrices that correspond to $r$ (negative reciprocation) and $+1$ (adding one) in the
following sense.  We  define the fraction of a vector $v$, $[v],$ by the formula 

$$[v]=[  \left(
\begin{array}{cc}
     a \\
     b
\end{array}
\right)]
= \frac{a}{b}.$$

\noindent We also define the two basic matrices 
$$M(r) =  \left(
\begin{array}{cc}
     0 & -1 \\
     1 & 0
\end{array}
\right) \ \mbox{and} \  M(+1) =  \left(
\begin{array}{cc}
     1 & 1 \\
     0 & 1
\end{array}
\right) .$$
\noindent  Then $$[M(r)\cdot v] = \frac{-1}{[v]}  \ \mbox{and} \  [M(+1)\cdot v] = [v] + 1$$ 
\noindent  for any  vector $v$. So, we showed here that
addition of $[+1]$ and inversion suffice for generating all rational tangles.  By the result of 
this section, the players of the Square Dancing can dance their way through the intricacies of
$SL(2,\ZZ).$ 

\bigbreak

\noindent  {\bf History of rational knots and rational tangles.} \ 
As explained in  
\cite{Gr}, rational knots and links were first considered by O. Simony \cite{Si1, Si2, Si3, Si4}
in 1882, taking twistings and knottings of a band. Simony \cite{Si2} is  the first one to relate
knots to continued fractions. After about sixty years Tietze wrote  a series of papers
\cite{Ti1, Ti2, Ti3, Ti4}  with reference to Simony's work. 
Reidemeister \cite{Rd2} in 1929 calculated the knot group of a special class of rational knots, 
but rational knots were studied  by Goeritz
\cite{Goe} and by Bankwitz and Schumann \cite{BS} in 1934. In \cite{Goe} and \cite{BS} rational 
knots are represented as plat closures of four-strand braids. 

Figure 2 in \cite{BS}  illustrates a rational
tangle, but no special importance is given to this object. The rational tangle is obtained by a
four-strand braid by plat-closing only the top four ends. A rational tangle obtained this way may 
be said to be between the twist form (Definition 1) and the standard form (Definition 4), in the
sense that, if we twist neighbouring endpoints starting from two trivial arcs, we may twist  to the 
right and to the left but only to the bottom, not to the top (see Lage 3 of \cite{BS}). In
\cite{Goe} and
\cite{BS}  proofs are given independently and with different techniques that rational knots have
$3$-strand-braid representations (in \cite{BS} using the horizontal-vertical structure of the 
rational tangles), in the sense that the first strand of the four-strand braids can be free of
crossings.  The
$3$-strand-braid representation of a four-plat corresponds to the numerator of a rational tangle in
standard form. In \cite{Goe} and \cite{BS}  proofs are also given that rational knots are 
alternating. The proof of this fact in \cite{BS} can be easily applied on the corresponding
rational tangles in standard form.  

It was not until 1956 that Schubert \cite{Sch} 
classified rational knots by finding canonical forms via representing them as $2$-bridge knots. 
His proof was based on Seifert's observation that the $2$-fold branched coverings of
$2$-bridge knots give rise to lens spaces and on the classification of lens spaces  by
Reidemeister  \cite{Rd3} using Reidemeister torsion (and later by Brody
\cite{Br1, Br2} using knot theory in lens spaces). See also \cite{PY}. Schubert's theorem was  reformulated
by Conway
\cite{C1} in  terms of  rational tangles.  See the paper
of Siebenmann \cite{Sie} for an excellent exposition and see the book by Bonahon and Siebenmann 
\cite{BSI}  for developments about tangles circa 1980. 
\bigbreak

\noindent {\bf Acknowledgments.}  It gives us great pleasure to thank John Conway, 
Ray Lickorish and Jozef Przytycki for useful
conversations. We also thank Jozef Przytycki for telling us about the work of Simony and of Tietze.


\bigbreak

\noindent {\sc L.H.Kauffman: Department of Mathematics, Statistics and Computer Science, University
of Illinois at Chicago, 851 South Morgan St., Chicago IL 60607-7045, U.S.A.}

\vspace{.1in}
\noindent {\sc S.Lambropoulou: National Technical University of Athens, Department of Mathematics,
Zografou campus, GR-157 80 Athens, Greece. }

\vspace{.1in}
\noindent {\sc E-mails:} \ {\tt kauffman@math.uic.edu  \ \ \ \ \ \ \ \ sofia@math.ntua.gr 

\noindent http://www.math.uic.edu/$\tilde{~}$kauffman/ \ \ \ \ http://users.ntua.gr/sofial}

\end{document}